\documentclass[review,onefignum,onetabnum]{siamart171218}

\usepackage{lipsum}
\usepackage{amsfonts}
\usepackage{graphicx}
\usepackage{xcolor}
\usepackage{hyperref}
\usepackage{subfigure}

\usepackage{epstopdf}
\usepackage{algorithmic}

\usepackage{bm}
\usepackage{tikz}
\usetikzlibrary{decorations.pathreplacing,calligraphy}

\newcommand{\bU}{{\bm{U}}}
\newcommand{\bv}{{\bm{v}}}
\newcommand{\bxi}{\bm{\xi}}
\newcommand{\bd}{{\bm{d}}}

\newcommand{\bF}{{\bm{F}}}
\newcommand{\bE}{{\bm{E}}}

\usepackage{ntheorem}

\theoremstyle{plain}
\theoremheaderfont{\normalfont\bfseries}
\theorembodyfont{\normalfont}
\theoremseparator{.}
\newtheorem{example}{Example}[section]

\ifpdf
  \DeclareGraphicsExtensions{.eps,.pdf,.png,.jpg}
\else
  \DeclareGraphicsExtensions{.eps}
\fi

\newsiamremark{remark}{Remark}
\newsiamremark{hypothesis}{Hypothesis}
\crefname{hypothesis}{Hypothesis}{Hypotheses}
\newsiamthm{claim}{Claim}

\headers{Conservative SL FD Scheme Using GNN}{Y. Chen, W. Guo, and X. Zhong}

\title{Conservative Semi-Lagrangian Finite Difference Scheme for Transport Simulations Using Graph Neural Networks\thanks{Submitted to a SIAM Journal.
\funding{W. Guo is partially supported by the NSF grant NSF-DMS-2111383, Air Force Office of Scientific Research FA9550-18-1-0257.  X. Zhong is partially supported by the NSFC Grant 12272347. }}}

\author{Yongsheng Chen\thanks{School of Mathematical Sciences, Zhejiang University, Hangzhou, 310027, China.
  (\email{22035024@zju.edu.cn}).}
\and 
Wei Guo
\thanks{Corresponding author. Department of Mathematics and Statistics, Texas Tech University, Lubbock, TX, 70409, USA. 
(\email{weimath.guo@ttu.edu}). }
'\and
Xinghui Zhong
\thanks{ School of Mathematical Sciences, Zhejiang University, Hangzhou, 310027, China. (\email{ zhongxh@zju.edu.cn})}
}
\usepackage{amsopn}

\makeatletter
\newcommand*{\addFileDependency}[1]{% argument=file name and extension
  \typeout{(#1)}% latexmk will find this if $recorder=0 (however, in that case, it will ignore #1 if it is a .aux or .pdf file etc and it exists! if it doesn't exist, it will appear in the list of dependents regardless)
  \@addtofilelist{#1}% if you want it to appear in \listfiles, not really necessary and latexmk doesn't use this
  \IfFileExists{#1}{}{\typeout{No file #1.}}% latexmk will find this message if #1 doesn't exist (yet)
}
\makeatother

\newcommand*{\myexternaldocument}[1]{%
    \externaldocument{#1}%
    \addFileDependency{#1.tex}%
    \addFileDependency{#1.aux}%
}
\ifpdf
\hypersetup{
  pdftitle={A graph-based conservative semi-Lagrangian finite difference scheme for transport simulations},
  pdfauthor={Y. Chen, W. Guo, X. Zhong}
}
\fi

\myexternaldocument{ex_supplement}

\begin{document}
\nolinenumbers

\maketitle

% REQUIRED
\begin{abstract}
Semi-Lagrangian (SL) schemes are highly efficient for simulating transport equations and are widely used across various applications. Despite their success, designing genuinely multi-dimensional and conservative SL schemes remains a significant challenge. Building on our previous work [Chen et al., J. Comput. Phys., V490 112329, (2023)], we introduce a conservative machine-learning-based SL finite difference (FD) method that allows for extra-large time step evolution. At the core of our approach is a novel dynamical graph neural network designed to handle the complexities associated with tracking accurately upstream points along characteristics. This proposed neural transport solver learns the conservative SL FD discretization directly from data,  improving accuracy and efficiency compared to traditional numerical schemes, while significantly simplifying algorithm implementation. We validate the method's effectiveness and efficiency through numerical tests on benchmark transport equations in both one and two dimensions, as well as the nonlinear Vlasov-Poisson system.

\end{abstract}

% REQUIRED

% REQUIRED

\section{Introduction}
Transport equations are prevalent in various scientific and engineering disciplines, including numerical weather prediction \cite{lynch2006weather,bauer2015quiet}, climate modeling \cite{schneider2017climate,neumann2019assessing}, and plasma physics \cite{tang2005advances}, among many others. The last several decades have witnessed tremendous development of effective computational tools for simulating transport equations, while several challenges remain in approximating nonsmooth or multi-scale structures with high order accuracy and robustness, in simulating large-scale problems efficiently with manageable resources, and in preserving inherent structures of the underlying equations.

Semi-Lagrangian (SL) schemes \cite{wiin1959application,staniforth1991semi,sonnendrucker1999semi} are a popular numerical tool for simulating transport equations, offering numerous computational advantages. As a mesh-based approach,  SL schemes are able to support various spatial discretization framework, including   finite difference (FD) methods \cite{carrillo2007nonoscillatory,qiu2010conservative,qiu2011conservative,lentine2011unconditionally}, finite volume (FV) methods \cite{filbet2001conservative,lauritzen2010conservative,zheng2022fourth}, and discontinuous Galerkin finite element methods \cite{qiu2011positivity,rossmanith2011positivity,guo2014conservative,cai2017high,einkemmer2020semi}. Furthermore, SL schemes evolve grid-based numerical solutions by following the characteristics, enjoying the beneficial properties of both Eulerian and Lagrangian approaches. Hence, such schemes may avoid the statistical noise and  time step restrictions simultaneously, resulting in significant efficiency. Moreover, due to their distinctive unconditional stability, SL schemes are capable of conveniently bridging the disparate time scales in the problem \cite{liu2021self}, which is highly desirable for multi-scale transport simulations. However, designing genuinely multi-dimensional and conservative SL schemes remains a significant challenge. Hence, a splitting approach is often incorporated to circumvent the difficulty at the cost of introducing splitting errors \cite{cheng1976integration}. Meanwhile, similar to other conventional grid-based schemes, the accuracy of SL solvers is typically constrained by the resolution of the simulation mesh. Therefore, high-resolution grids are necessary to accurately capture the fine-scale solution structures of interest, resulting in significant computational costs, especially for large-scale simulations.

With the rapid development in machine learning (ML) and computing power over recent decades, the integration of ML tools with the simulation of partial differential equations (PDEs) has emerged as a thriving field. Various ML-based PDE solvers are developed to address inherent shortcomings of  traditional numerical schemes, leveraging the expressive power of neural networks (NNs) and advancements in automatic differentiation technology \cite{Ilya_fix_2017}. These neural PDE solvers have achieved tremendous success across numerous applications, demonstrating improved efficiency and accuracy compared to traditional solvers. Among these developments, one notable example is the physics informed neural networks (PINNs) \cite{raissi_physics_2017,raissi_physics_2017-1}, where the solutions are parameterized with  NNs and trained using a physics-informed loss function. PINNs are extensively employed to solve complex problems in various fields \cite{lu_deepxde_2021,lu_physics-informed_2021,mcclenny_self-adaptive_2022,pang_fpinns_2019,raissi_physics-informed_2019,yu2022gradient,jagtap_conservative_2020}, and a comprehensive review of the literature on PINNs was provided in \cite{cuomo2022scientific}. Another category of NN-based PDE solvers, known as neural operators, focuses on learning the mapping from initial conditions to solutions at  later time 
$t$. Some related works include \cite{lu_deeponet_2021,Li_fourier_2021,kissas2022learning,bhattacharya_model_2021,li_neural_2020}. In addition, autoregressive methods \cite{bar2019learning,hsieh2019learning,brandstetter2021message,chen_cell-average_2022} offer a different approach specifically designed for time-dependent problems. These methods simulate the PDEs iteratively, resembling conventional numerical methods that employ time marching.

Recently, graph neural networks (GNNs) have gained significant research attention as a powerful approach for simulating complex physical systems due to their superior learning capabilities, flexibility, and generalizability. They excel particularly in adapting to unstructured grids and high-dimensional problems. In particular, existing GNN-based PDE solvers can be roughly classified into two groups. The first group focuses on learning the continuous graph representation of the computational domain. Related works include utilizing GNNs to learn a continuous mapping \cite{li_neural_2020,you2022nonlocal}, learning the solutions of different resolutions \cite{li2020multipole}, developing a continuous time differential model for the dynamical systems \cite{iakovlev2020learning}, and integrating the GNN with differential PDE solvers \cite{belbute2020combining}. In such representation, the architecture is invariant to data resolution, and the interactions between nodes are enhanced, which is desired for simulating complex physical systems. The second group  explores the capabilities of GNNs in learning mesh-based simulations. A notable example is the MeshGraphNet \cite{pfaff2020learning,fortunato2022multiscale}, which encodes the mesh information and the corresponding physical parameters into a graph. Other GNN related works include the particle-based approach \cite{sanchez2020learning}, novel GNN architectures \cite{gladstone2023gnn} designed to facilitate long-range information exchange,  and a message-passing framework \cite{brandstetter2021message} for  PDE simulations. These graph-based approaches can accurately simulate complex physical systems and exhibit remarkable generalization properties.

  % A particle-based approach was adopted to simulate a range of physics phenomena \cite{sanchez2020learning}. Meanwhile, to facilitate the distant exchange of information across the computational domain, two novel GNN architectures were proposed \cite{gladstone2023gnn} to overcome this challenge. The consideration of boundary conditions saw the construction of a physics-embedded neural network to predict long-term behaviors \cite{horie2022physics}. Furthermore, a comprehensive framework based on message passing was proposed \cite{brandstetter2021message} to solve PDEs. These graph-based approaches can accurately simulate physics systems and exhibit impressive generalization properties.

In this paper, we develop an ML-based conservative SL FD scheme which improves  traditional SL schemes as well as existing ML-based transport solvers. Our method belongs to the category of autoregressive methods and aims to learn the optimal SL discretization through a data-driven approach. To achieve this, we propose an end-to-end neural PDE solver that consists of three classic network blocks: the encoder, the processor, and the decoder, inspired by the works developed in  \cite{sanchez2020learning,battaglia2018relational,brandstetter2021message,equer2023multi}. The encoder, constructed as a convolutional neural network (CNN) \cite{lecun1995convolutional,lecun1998gradient},  is designed to learn the embedding of node features. It incorporates the normalized shifts as part of the input of the NN, observing the fact that such quantities contribute to computing the traditional SL FD discretization but in a rather complex manner. For the processor, a \emph{dynamical} graph is generated based on the normalized shifts, allowing us to effectively handle the geometries associated with upstream points tracking. Subsequently, we construct a GNN which is composed of several graph convolutional layers to process and propagates the information. Through message passing, each node can yield a latent feature vector. The decoder, consisting of a multilayer perceptron (MLP) and a constraint layer, utilizes these feature vectors to predict the SL discretization and  guarantee exact mass conservation. Note that we  replace the most expensive and complex component of the SL formulation with an end-to-end data-driven approach. This eliminates the need for explicit implementation of tracing upstream points, resulting in improved efficiency and significantly reducing the human effort required for implementation. Meanwhile, the proposed dynamical graph approach enables a local interpolation procedure, thereby allowing for large time steps for evolution, in contrast to our previous work \cite{chen2023learned}.  As with other ML-based discretization methods \cite{bar2019learning,zhuang2021learned,kochkov_machine_2021,chen2023learned}, the learned SL discretization with a coarse grid can accurately capture fine-scale features of the solution, which often demands much finer grid resolution for a traditional polynomial-based discretization, leading to significant computational savings. %The use of CNN and GNN can also incorporate inductive bias and give the model improved generalization.
Furthermore, we extend this method for solving the nonlinear Vlasov-Poisson (VP) system. The inherent nonlinearity of the VP system presents an additional challenge in accurately tracking the underlying characteristics. To overcome this difficulty, we integrate the high-order Runge-Kutta (RK) exponential integrators (RKEIs), introduced in \cite{celledoni2003commutator,cai2021high}. By employing the RKEI, the VP system can be decomposed into a sequence of linearized transport equations, each with a constant frozen velocity field \cite{cai2021high,zheng2022fourth}. This decomposition allows us to apply the proposed GNN-based SL scheme, resulting in a novel data-driven conservative SL FD VP solver without operator splitting.

% The outline is not required, but we show an example here.
 The rest of the paper is organized as follows. In Section \ref{sec:algorithm}, we provide background and brief reviews of related works, including a traditional SL FD scheme for one-dimensional (1D) transport equations, along with recent developments in neural solvers for time dependent PDEs. In Section \ref{sec:main}, we introduce the proposed GNN-based conservative SL FD scheme with application to the VP system. In Section \ref{sec:experiments}, the numerical results are presented to demonstrate the performance of the proposed method. The conclusion and future works are discussed in Section \ref{sec:conclusions}.

\section{Background and related work}
\label{sec:algorithm}
\subsection{Semi-Lagrangian finite difference scheme}\label{subsec:slfv}
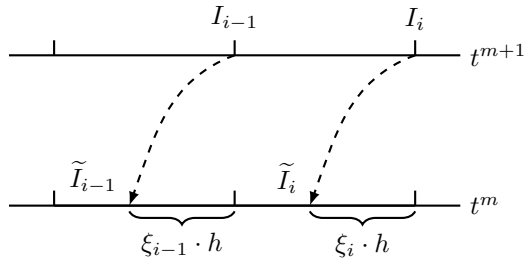
\begin{figure}[h!]
	\centering
		\begin{tikzpicture}[x=1cm,y=1cm]%[scale=1.0]
		\begin{scope}[thick]
		% \draw (-3,3) node[fill=white] {};
		% \draw (-3,-1) node[fill=white] {};
		%\draw (60:-1cm) node[fill=white] {$E$} -- (60:3cm) node[fill=white] {$F$};
		\draw[black]  (-3,0) -- (3,0)
		node[right]{$t^{m}$};
		\draw[black] (-3,2) -- (3,2)
		node[right]{$t^{m+1}$};
		
		%\draw[snake=ticks,segment length=2.2cm] (-2.4,0) -- (2.8,0);
		% \draw[snake=ticks,segment length=2.4cm] (-2.4,2) -- (0,2) node[above] {$x_{i-\frac12}$};
		% \draw[snake=ticks,segment length=2.4cm] (0,2) -- (2.4,2) node[above] {$x_{i+\frac12}$};
		
		\draw[black] (-2.4,0) -- (0,0);
		\draw[black] (0,0) -- (2.4,0);
		\draw[black] (0,2) -- (0,2.2)node[above] {$I_{i-1}$};
            \draw[black] (0,0) -- (0,0.2);
		\draw[black] (2.4,2) -- (2.4,2.2)node[above] {$I_{i}$};
            \draw[black] (2.4,0) -- (2.4,0.2);
		\draw[black] (-2.4,2) -- (-2.4,2.2);
            \draw[black] (-2.4,0) -- (-2.4,0.2);            
		\draw[-latex,dashed](0,2) to[out=200,in=70] (-1.4,0) node[above left=1pt] {$\widetilde{I}_{i-1}$};  
		\draw[-latex,dashed](2.4,2) to[out=200,in=70] ( 1.,0) node[above left=1pt] {$\widetilde{I}_{i}$};
            % \draw[decorate, decoration = {brace}](-1.4,0) -- (0,0) 
		% \draw[snake=brace,mirror snake,red,thick] (-1.4,0) -- (0,0) node[below left=2.5pt] {$\widetilde{I}_{i,i-1}$};
		% \draw[snake=brace,mirror snake,red,thick] (0,0) -- (1,0) node[below left=2.5pt] {$\widetilde{I}_{i,i}$};
  \draw [decorate,decoration={brace,amplitude=5pt,mirror,raise=0.5ex}]
  (-1.4,0) -- (0,0)node[midway,yshift=-1.5em] {$\xi_{i-1}\cdot h$};
  \draw [decorate,decoration={brace,amplitude=5pt,mirror,raise=0.5ex}]
  (1.,0) -- (2.4,0)node[midway,yshift=-1.5em] {$\xi_{i}\cdot h$};		
		\end{scope}
		\end{tikzpicture}

	\caption{Schematic illustration of the 1D SL FD scheme.\label{fig:slfv}}	
\end{figure}

In this section, we review the  SL FD scheme  for linear transport equations proposed in \cite{qiu2011conservative}. We start with the following 1D equation in the conservative form 
\begin{equation}
    \label{eq:transport1d}
    u_t + (a(x,t)u)_x = 0,\quad x\in\Omega,
\end{equation}
where $a(x,t)$ is the velocity function. For simplicity, consider a uniform partition of the domain $\Omega$ with $N$ grid points, denoted as $\{I_i\}_{i=1}^N$, where each grid point $I_i$ has a coordinate $x_i$. Denote the mesh size as $h=x_i-x_{i-1}$.  The numerical solution  $\{U_i^m\}_{i=1}^N$ approximates the solution value $u(x_i,t^m)$ at node $I_i$  and  time step $t^m$. The information of the equation \eqref{eq:transport1d} propagates according to the characteristic equation
\begin{equation}
\label{eq:characteristic}
    \frac{dx(t)}{dt} = a(x(t),t).
\end{equation}
To update the solution in the SL setting, we evolve \eqref{eq:characteristic} backwards from $t^{m+1}$ to $t^{m}$ at each grid point $I_i$, and obtain the corresponding upstream point $\widetilde{I}_i$ with the coordinate $\widetilde{x}_i$, as shown in Figure \ref{fig:slfv}. Define the normalized shift as
\begin{equation}
\label{eq:shift}
\xi_i = \frac{\widetilde{x}_i - x_i}{h},\quad i=1,\ldots,N,
\end{equation}
which plays a key role in the algorithm development.
%
%
% The SL FD method updates the  point value approximation  to next time step $t^{m+1}$ by tracking the characteristics governed by the ordinary differential equation
% % \begin{equation}
% % \label{eq:characteristic}
% %     \frac{dx(t)}{dt} = a(x(t),t).
% % \end{equation}
% Compared to SL FV schemes, the design of a conservative   SL FD scheme is more involved. 

%There is much less research and development of conservative SL FD schemes  compared to the SL FV counterpart in the literature.  
While  conservative SL FV schemes are well-developed in the literature, the SL FD counterparts received considerably less research attention. In the seminal work \cite{qiu2011conservative},  a flux difference SL formulation was proposed
\begin{equation}
\label{eq:slfd1d}
U_i^{m+1} = U_i^{m} - \frac{1}{h} (\hat{F}_{i+\frac12} - \hat{F}_{i-\frac12}),\quad i=1,\ldots,N,
\end{equation} 
where the numerical flux $\hat{F}_{i+\frac12}$ relies on a two-step reconstruction procedure. Interested readers are referred to \cite{qiu2011conservative} for a comprehensive algorithmic description. Moreover, the desired mass conservation property is automatically attained due to the flux difference form. 
%
% where the flux difference $\frac{1}{h} (\hat{F}_{i+\frac12} - \hat{F}_{i-\frac12})$ approximates $\int_{t^m}^{t^{m+1}}(a(x_i,t)u(x_i,t))_xdt$. 
Note that $\hat{F}_{i+\frac12}$ is fully determined by the reconstruction stencil, the normalized shifts $\{\xi_i^m\}_{i=1}^N$, and the point values $\{U_i^m\}_{i=1}^N$. For example, assume $-1\le\xi_{i-1},\xi_i\le 0$, the numerical fluxes can be approximated by
$$
\hat{F}_{i-\frac12} = -h\xi_{i-1}U_{i-1}^m,\quad \hat{F}_{i+\frac12} = -h\xi_iU_i^m
$$
with first order accuracy. 
Then, the solution is updated by
\begin{equation}
\label{eq:1st}
U_i^{m+1} =  U_i^m - (-\xi_i U^m_i + \xi_{i-1} U^m_{i-1}) =-\xi_{i-1} U^m_{i-1} +  (1+\xi_i)U^m_i,
\end{equation}
using the stencil $\{I_{i-1},I_{i}\}$. It is important to highlight that such an SL FD formulation allows for large time step evolution,  resulting in considerable computational efficiency. This is achieved by utilizing the stencils centered around the positions of the upstream points to update solutions at grid points, which is known as \emph{local} reconstruction or interpolation.  In addition, the high order accuracy can be achieved by employing a wider stencil for reconstruction. Similar to the first order scheme \eqref{eq:1st}, we can express a high order SL FD scheme as follows 
\begin{equation}
\label{eq:slfd}
U_i^{m+1} = \sum_{j\in\mathcal{N}_{in}(i)} d^m_{ji}U_j^m,
\end{equation}
where $\mathcal{N}_{in}(i)$ denotes the stencil used to update the solution $U_i^{m+1}$. If the coefficients $\{d^m_{ji}\}$ depend solely on the normalized shifts as well as the used stencils,  without  reliance on the numerical approximations $\{U_i^m\}$, then the scheme \eqref{eq:slfd} is  linear. It is possible to incorporate the nonlinear WENO mechanism in the reconstruction \cite{jiang1996efficient}, as outlined in \cite{qiu2011conservative}, where the coefficients $\{d^m_{ji}\}$ also depend on  $\{U_i^m\}$. The resulting SL FD WENO scheme enjoys improved stability and robustness when approximating  discontinuous solutions.  

It is shown  in  \cite{chen2023learned} that an SL scheme expressed in the form of \eqref{eq:slfd}  for solving \eqref{eq:transport1d} is mass conservative, i.e., $\sum_{i=1}^N U_i^{m+1} = \sum_{i=1}^N U_i^{m}$, if
\begin{equation}\label{eq:mass}
    \sum_{i\in \mathcal{N}_{out}(j)}d^m_{ji}=1,\quad \forall j,
\end{equation}
where $\mathcal{N}_{out}(j)$ denotes a collection of grip points $\{I_i\}$ for which $U_j^m$  contributes to the update of $\{U_i^{m+1}\}$, i.e., the region of influence of $U_j^m$. Furthermore, if the scheme is linear, then \eqref{eq:mass} is also a necessary condition. Interestingly, it can be verified that the SL FD WENO scheme proposed in \cite{qiu2011conservative}, though conservative due to its flux difference form, does not satisfy the condition \eqref{eq:mass}. Nevertheless, \eqref{eq:mass} provides a means to ensure mass conservation for a transport scheme.  It was utilized in \cite{lentine2011unconditionally} to design a first order conservative SL FD scheme. 
%In instance, in \cite{lentine2011unconditionally},  \eqref{eq:mass} is utilized to design a first order conservative SL FD scheme. 

A significant limitation of the SL FD formulation is that its genuine extension to multi-dimensional equations is highly challenging. To circumvent this difficulty, one can perform dimensional splitting which decouples the underlying multi-dimensional equation into a series of 1D equations, with the trade-off of incurring splitting errors \cite{cheng1976integration,qiu2011conservative}. Furthermore, to the authors' best knowledge, there does not exist a conservative non-splitting multi-dimensional SL FD scheme with accuracy  higher than first order in the literature.
On the other hand, the SL FV architecture is capable of overcoming such a limitation and being non-splitting, conservative, and unconditionally stable simultaneously, see, e.g., \cite{lauritzen2010conservative}. However, this approach necessitates accurate tracking of deformed upstream cells, which is  notably intricate and demands substantial human effort to implement. Moreover, extending such a non-splitting SL FV method to three-dimensional or higher dimensional problems is not directly achievable. Below, we will leverage cutting-edge ML techniques and develop a genuine multi-dimensional conservative SL FD scheme.

% The upstream point $\widetilde{I}_i$ of grid point $I_i$ is defined by evolving \eqref{eq:characteristic} with final values $x(t^{m+1}) = x_{i}$ backward to $t^m$, as shown in Figure \ref{fig:slfv}, and the solution is the coordinate of $\widetilde{I}_i$ denoted by  $\widetilde{x}_i$.

\subsection{Autoregressive neural transport solvers}

Recently, the rapid development of deep NNs and their superior ability to approximate complex  functions  have inspired researchers to design innovative numerical methods for solving PDEs, often referred as neural solvers in the literature. These approaches leverage the data (e.g., from classical numerical solvers or experiments), to train the underlying neural solver, aiming for improved accuracy and efficiency when generalized  to new context. In particular, for time-dependent problems, an autoregressive neural method iteratively updates the approximate solution in a manner similar to conventional numerical methods with time marching, see \cite{brandstetter2021message}. By contrast, neural operator methods aim to learn the mapping from initial conditions to solutions at later time, see, e.g \cite{lu_deeponet_2021,Li_fourier_2021,kissas2022learning}.

A noteworthy development in autoregressive neural methods is the ML-based discretization approach, see \cite{pathak2020using,obiols2020cfdnet,sirignano2020dpm,tompson2017accelerating,um2020solver}. This approach replaces the polynomial-based interpolation/reconstruction  typically employed by traditional numerical methods with approximations by deep NNs. By training on high-quality data, the NN learns the optimal numerical discretization specific for the underlying PDE, leading to significant computational efficiency. For time-dependent PDEs, the ML-based discretization updates the solution approximations iteratively using the method of lines. Moreover, the ML-based discretization 
inherits its structure from the classical solver. Hence, it facilitates the convenient enforcement of inherent physical constraints, such as  conservation of mass, momentum, and energy, at the discrete level, which is crucial for generalization and reliability of the neural solver \cite{zhuang2021learned,kochkov_machine_2021,bar2019learning}.  Recently, we developed a conservative ML-based SL FV scheme, which is able to take larger time steps for evolution by incorporating information of characteristics as inductive bias \cite{chen2023learned}. Meanwhile, this method employs a set of \emph{fixed} stencils to update the solution, and thus is constrained by a CFL condition.

Another group of autoregressive neural solvers is the message-passing neural PDE methods \cite{brandstetter2021message,equer2023multi}. Such an approach is based on message-passing GNN architecture, which offers considerable flexibility, and  employs the Encode-Process-Decode framework developed in \cite{sanchez2020learning,battaglia2018relational}. In particular, the method models the mesh as a graph. To update the solution, it begins by encoding the state and PDE parameters into the graph structure using an encoder. During the message-passing phase, it incorporates the relative positions of nodes together with the solution differences, allowing the model to exploit the translational symmetry of the underlying solution structures. This process also represents an analogy to traditional numerical differential operators.
Finally, it employs a decoder to extract the node information and update the solutions based on the learned representations. Such a message-passing neural solver, employing GNN, presents an effective end-to-end data-driven framework for simulating PDEs and enjoys enhanced flexibility for generalization. However, unlike the ML-based discretization approach, this method does not explicitly account for physical constraints, such as mass conservation.   In addition, for solving transport equations, the time step is also restricted by an CFL condition.
%it is also limited by a CFL-type time step restriction. 

In this study, we focus on first-order transport equations and introduce a novel autoregressive neural solver. Specifically, we aim for an ML-based discretization method that achieves higher efficiency than existing autoregressive neural solvers.  The proposed method adheres to the fundamental principle of transport equations, namely propagating information  along characteristics. In particular, it incorporates the SL mechanism to update solutions, building upon our previous work \cite{chen2023learned}. Furthermore, we develop a novel Encode-Process-Decode architecture and employ a message-passing GNN to handle the complexities associated with tracking upstream points. In addition to the computational efficacy by  existing autoregressive neural solvers, our approach enables large time step evolution while ensuring exact local mass conservation. In \cite{larios2022error}, an ML-based SL approach was developed in the context of the level-set method, which employs a NN to correct the local error incurred by the standard SL FD method. While this approach is able to avoid the CFL time step restriction, it does not conserve mass.

\section{GNN-based conservative SL FD scheme}\label{sec:main}
In this section, we formulate a novel GNN-based SL FD scheme, which achieves exact mass conservation and unconditional stability simultaneously. As mentioned above, for an SL scheme to achieve unconditional stability, the underlying interpolation/reconstruction process must be localized concerning the positions of upstream points. The primary challenge for designing the algorithm stems from potential irregular geometry and distant positioning of upstream points with respect to the original grid points, especially when employing significantly large time steps.
To tackle these complexities, we propose a Encode-Process-Decode framework. This framework centers around a  novel dynamical GNN processor that captures the geometric relationships between the background and upstream grid points, thus effectively  facilitating the propagation of information following characteristics.

\subsection{Architecture}\label{sec:alg-linear}
\begin{figure}[htbp]
  \centering
  \label{fig:network}\includegraphics[width=1\textwidth]{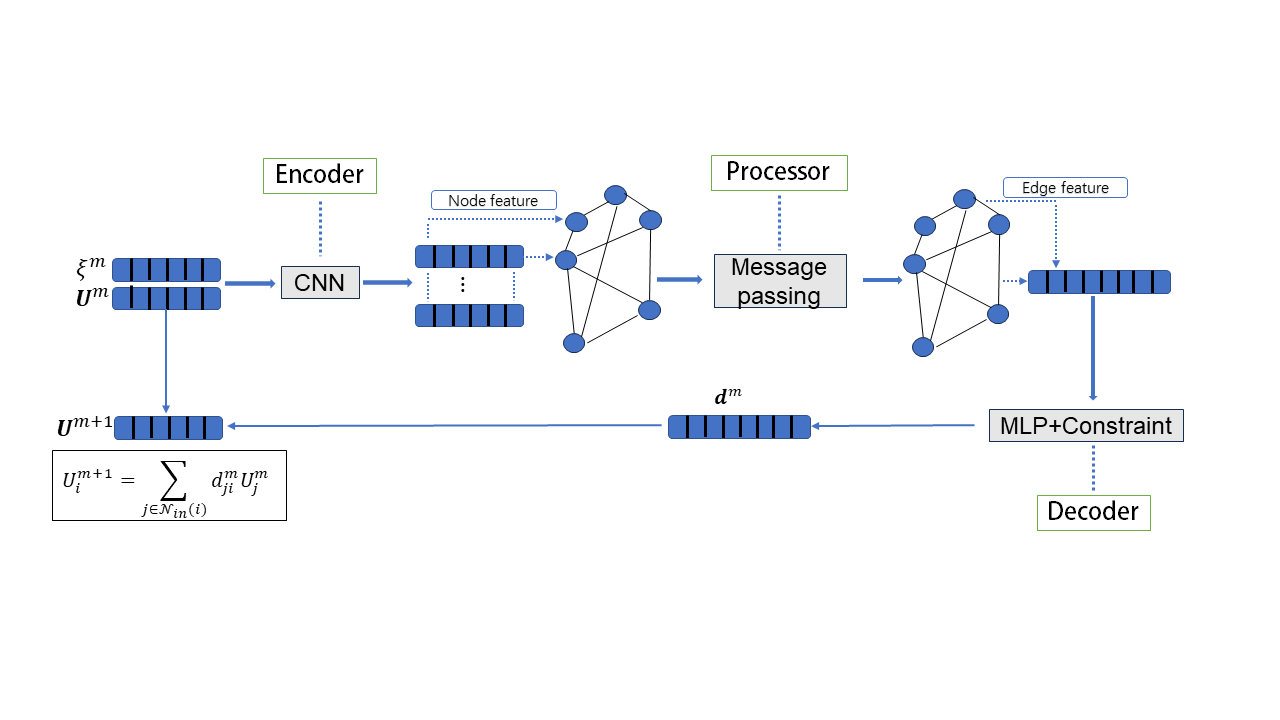}
   
  \caption{Illustration of the proposed GNN-based SL FD method.}
\end{figure}

For simplicity, we illustrate the main architecture of the proposed method for the 1D case and briefly discuss the extension to the two-dimensional (2D) case at the end of this subsection. We denote by $\bU^m$, $\bxi^m$, and $\bd^m$ the collections of $\{U^m_i\}$, $\{\xi_i^m\}$, and $\{d^m_{ji}\}$, respectively.  We adopt the effective Encode-Process-Decode framework \cite{sanchez2020learning,battaglia2018relational}, with crucial adjustments tailored for transport equations. The proposed GNN-based SL FD method is schematically illustrated in Figure \ref{fig:network}. %The input and output of each block and its functions are defined as follows.\\

\noindent{\bf Encoder.} The encoder, constructed as a CNN,  is employed to compute the embedding of node features with input  $\bU^m$ and $\bxi^m$. 
% Here the node feature $x_i = (U_i^m, \xi_i^m)\in \mathbb{R}^2$. 
It is proposed in our previous work \cite{chen2023learned}, which is motivated by the observation that the coefficients $\bd^m$ are fully determined by  $\bU^m$ and $\bxi^m$, see e.g., \eqref{eq:1st}. In particular, the encoder maps  node features $(\bU^m,\bxi^m)$  to   embedding vectors $ \bm{h}^{m,(0)}$:
\begin{equation*}
    \bm{h}^{m,(0)} = Enc(\bU^m,\bxi^m),
\end{equation*}
where the CNN encoder $Enc$ takes $\bU^m$ and  $\bxi^m$ as a two-channel input and is constructed as a stack of 1D convolutional layers and nonlinear activation functions, such as ReLU. The employed CNN can effectively capture hierarchical features of the solution \cite{lecun1989handwritten,lecun1995convolutional}, which are highly desirable for transport modelling. 
Note that other NN architectures for embedding, such as U-Net \cite{ronneberger2015u} and Vision Transformer \cite{han2022survey}, can also be utilized for improved performance. Furthermore, when working with an unstructured mesh, a GNN can be employed as the encoder to embed the node features, similar to the message-passing neural PDE solver \cite{brandstetter2021message}. \\

\noindent{\bf Processor.}  Although CNNs can efficiently extract features from solution structures, the flexibility is limited by  the fixed size of their convolution kernels.  In particular,  with CNN, only fixed stencils are allowed to design the SL discretization, resulting in the undesired CFL time step restriction, see \cite{chen2023learned}. Hence,  for large time step evolution, we must employ very wide stencils to encompass the domain of dependence, leading to increased computational cost. Compared to CNNs, a GNN offers significantly more flexibility by enabling direct information propagation between nodes through the creation of corresponding edges, regardless of their geometric locations. Thus, the GNN architecture is well-suited for efficiently processing  possible irregular geometric information among upstream points and grid points, enabling  local interpolation required for large time step evolution.

\begin{figure}[!htbp]
 \centerline{\includegraphics[width=0.8\textwidth]{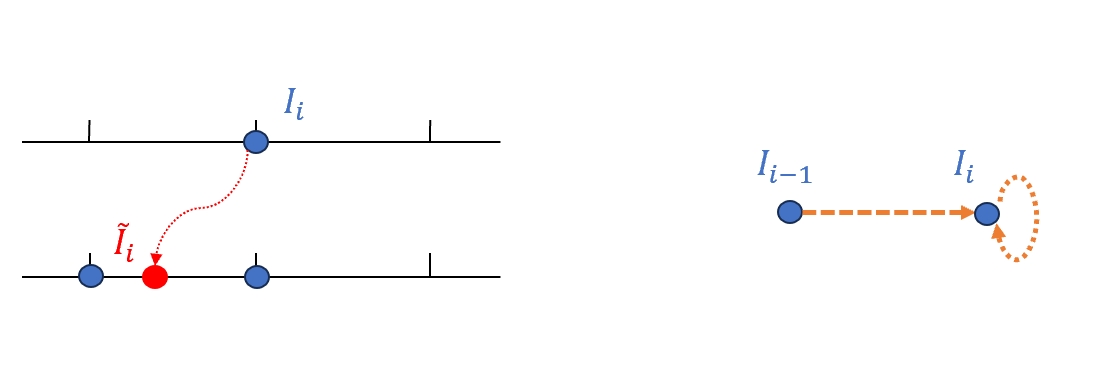}}
  \caption{Schematic illustration of the process of constructing the edges in 1D case. With local stencil $\mathcal{N}_{in}(i)=\{I_{i-1},I_i\}$ to update $U_i^{m+1}$, two edges $e_{i-1,i}$ and $e_{i,i}$ are created in the edge set $\mathcal{E}$.} \label{fig:edges_1d}
\end{figure}

In this paper, we develop a \emph{dynamical} graph architecture  $\mathcal{G}=(\mathcal{V},\mathcal{E})$,  with node embedding features $\bm{h}^{m,(0)}$ from the encoder, to capture the relationship between the upstream points and grid points. Here, a node $I_i \in \mathcal{V}$ represents a grid point, and an edge $e_{ji} \in \mathcal{E}$ represents a directed connection from node $I_j$ to node $I_i$. The graph construction process consists of two steps. First, for each node $I_i\in \mathcal{V}$,  the corresponding upstream point $\widetilde{I}_i$ is identified through characteristic tracing during the computation of $\xi_i^m$. Subsequently, we select a local stencil for interpolation. For instance, if the upstream point $\widetilde{I}_i$ falls in the interval $(x_{j-1},x_j]$, then a local two-point stencil $\{I_{j-1},I_j\}$ is chosen to update the solution, and we can create two directed edges $e_{j-1,i}$ and $e_{j,i}$ in the edge set $\mathcal{E}$.  Figure \ref{fig:edges_1d} illustrates a specific scenario where  the upstream point $\widetilde{I}_i$ is located within $(x_{i-1},x_i]$, and with the stencil $\{I_{i-1},I_i\}$, we create two edges $e_{i-1,i}$ and $e_{i,i}$. The latter, $e_{i,i}$, denotes a self-directed edge from node $I_i$ pointing to itself. Hence, unlike the case in \cite{brandstetter2021message}, the edges are not mesh edges. Instead, they represent the connections between upstream points and grid points in the stencil. Choosing a wider stencil can enhance the processor's generalization capability, but it increases costs due to the creation of more edges in the graph.
Let $\mathcal{N}_{in}(i)$ and $\mathcal{N}_{out}(i)$ denote the sets of nodes that have directed edges pointing directly towards and from the node $I_i$, respectively. Let $\mathcal{N}(i)=\mathcal{N}_{in}(i)\bigcup \mathcal{N}_{out}(i)$ denote the neighboring nodes of the node $I_i$. 

\begin{figure}[!htbp]
%\centering
 \centerline{\includegraphics[width=0.35\textwidth]{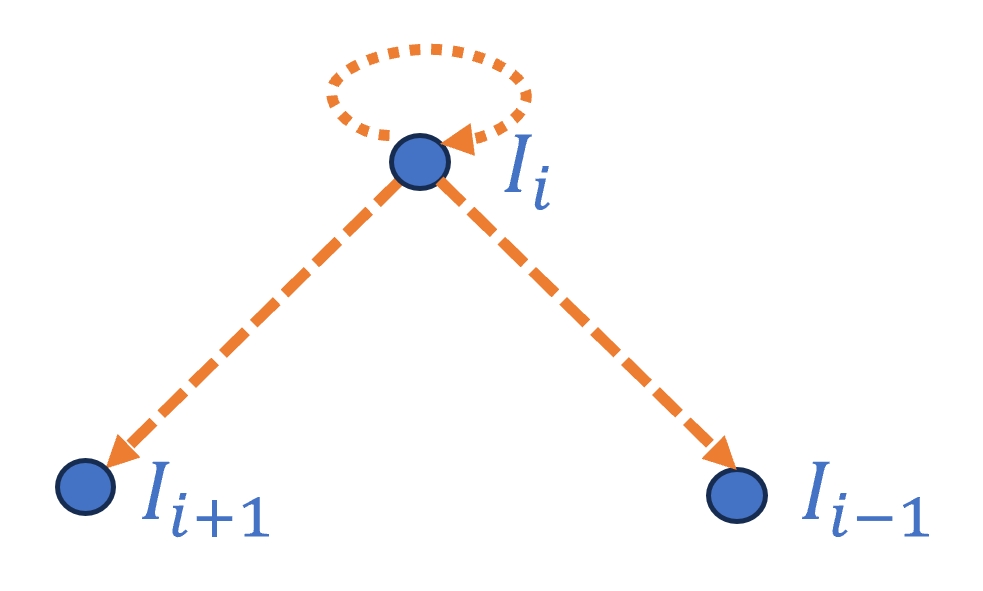}}
  \caption{Schematic illustration of $\mathcal{N}_{out}(i)=\{I_{i-1},I_i,I_{i+1}\}$, indicating 
 that $U_i^m$ will contribute to the computation of $U_{i-1}^{m+1}$, $U_{i}^{m+1}$, and $U_{i+1}^{m+1}$. For mass conservation, we need $d_{i,i-1}+d_{i,i}+d_{i,i+1}=1$.}\label{fig:edges_out}
\end{figure}

Noteworthy, for the variable-coefficient or nonlinear transport equations, the locations of the upstream points may vary over time. Consequently, the corresponding graph $\mathcal{G}$ must be reconstructed at every time step, making it dynamical. In addition, the size of $\mathcal{N}_{in}(i)$ is identical to the size of the interpolation stencil and remains fixed. For instance, in the case of Figure  \ref{fig:edges_1d}, the size is two. On the other hand, the size of $\mathcal{N}_{out}(i)$ may vary among nodes and  change over time. In Figure \ref{fig:edges_out}, we present an example of $\mathcal{N}_{out}(i)$.

% The node feature of node $i$ is $h^{m,(0)}_i$, the output from encoder. Through the process of solving $\bxi^m$ in Section \ref{subsec:slfv}, we can construct each edge on the graph. For each node $i$($I_i$) on the graph, we find the upstream point $\widetilde{I}_i$.
% Then we determine the two closest grid points $I_{i_1}$ and $I_{i_2}$ to $\widetilde{I}_i$ and build the edges from nodes $i_1,i_2$ to node $i$, which means that the solutions $U^m_{i_1}$ of node $i_1$ and $U^m_{i_2}$ of node $i_2$ are used to update the solution $U^{m+1}_i$ of node $i$ at the next time step. For example, in Figure \ref{fig:edges_1d}, we can construct edges $e_{i-1,i}$ and $e_{i,i}$ for the node i. Denote $\mathcal{N}_i$ as the neighborhood of node $i$ in the graph, we can obtain the neighborhood $\mathcal{N}_i$ of each node $i$ in this way. Note that the process of constructing the graph depends on $\bxi^m$, which in turn depends on the time step $t_m$ and the solution trajectory $\bU^m$, meaning that we are constructing a dynamic graph. The graph structure of different solution trajectories is different and even the same solution trajectory will show different structural characteristics at different time steps.

Once the directed graph $\mathcal{G}$ is constructed, a GNN processor is employed to handle and propagate the information. In particular, a GNN exploits the graph structure to learn expressive node features. It iteratively updates a node's feature by aggregating and transforming messages from its neighbors. The utilized GNN processor is composed of $K$ graph convolutional layers, with each layer executing one step of message passing. Each step is defined as follows: 
\begin{subequations}
\begin{align}
    f_{ij}^{m,(k+1)} & = \phi(h_i^{m,(k)},h_j^{m,(k)}),\quad e_{ij}\in \mathcal{E}, \\
    h_i^{m,(k+1)} & = \varphi\left(h_i^{m,(k)},\rho\left(\left\{f^{m,(k+1)}_{ij},j\in \mathcal{N}(i)\right\}\right)\right),
\end{align}
\end{subequations}
for $k=0,1,\dots,K-1$, where $\phi$ is the message function, $\rho$ is the aggregation function, and $\psi$ is the update function.
% \begin{equation}
%     f_{i,j}^{m,(k+1)}  = \phi(h_i^{m,(k)},h_j^{m,(k)}),e_{ij}\in \mathcal{E},
% \end{equation}
% \begin{equation}
% h_i^{m,(k+1)}  = \varphi(h_i^{m,(k)},\rho(\left\{f^{m,(k+1)}_{ij},j\in \mathcal{N}(i)\right\})), k=1,2,\dots,K
% \end{equation}
Here, we adopt graph attention network (GAT) \cite{velivckovic2018graph,brody2021attentive} to perform the updates, while other GNN architectures with various designs for message-passing functions are also applicable, such as  graph convolutional network  \cite{chen2020simple,wu2019simplifying}, graph isomorphism network \cite{xu2018powerful}. For GAT, the message function $\phi$ is an MLP, which generates  messages on an edge by combining  the features of the nodes it connects. The aggregation function $\rho$ is an attention layer that aggregates the message received by a node, taking into account the relative importance of each message  determined by the attention mechanism. The update function $\varphi$ is a linear transformation. It is worth emphasizing that the processor does not incorporate edge features. 
While we recreate the graph at each time step, we only need to update the edge set $\mathcal{E}$. All other components, including NNs $\phi$, $\rho$, and $\varphi$, remain unchanged. Therefore, the associated cost is negligible.

We can consider augmenting the latent node  feature $h_i^{m,(0)}$ with the local coordinate of the upstream point $I_i$, relative to the underlying stencil, as an additional inductive bias, inspired by  the traditional SL FD methodology, e.g., \cite{qiu2011conservative}. However, numerical evidence suggests that such an inductive bias does not qualitatively alter the performance of the proposed scheme. 
\\

\noindent{\bf Decoder.} After message passing, we obtain the latent feature vector $h_i^{m,(K)}$ for each node $I_i$. Then we concatenate the feature vectors of two adjacent nodes connected by an edge, and obtain the feature vector $w^{m}_{ij}$ for edge $e_{ij}$ as
\begin{equation*}
    w^{m}_{ij} = [h_i^{m,(K)};h_j^{m,(K)}],\quad e_{ij} \in \mathcal{E}.
\end{equation*}
 % We present a schematic diagram of this process in Figure \ref{fig:concat}. 
 The decoder $Dec$ consists of two components. The first component  is an MLP $g$ that takes the edge features $w^{m}_{ij}$ as input and compute the pre-processed coefficients $\widetilde{d}^m_{ij}$: $\widetilde{d}^m_{ij} = g(w^m_{ij})$. Denote by $\widetilde{\bd}^m = \{\widetilde{d}^m_{ij}|e_{ij}\in \mathcal{E}\}$ and $\bm{w}^m = \{w^m_{ij}|e_{ij}\in \mathcal{E}\}$. 
 With a slight abuse of notations, we write $\widetilde{\bd}^m = g(\bm{w}^m)$, meaning that the MLP $g$ is applied to each component of $\bm{w}^m$.

 % $\bd^m = \{d^m_i\}$, where $d^m_i=\{d^m_{ij},e_{ij}\in \mathcal{E}\}$ represents the features of edges starting from node $i$. ????? i to j? or j to i?

 To ensure mass conservation at the discrete level, which is critical for accurate and stable long term transport simulations, we further propose to post-process the coefficients  $\widetilde{\bd}^m$ and make sure that the sum of contributions of node $I_i$ that are employed to update the solutions is 1, see \eqref{eq:mass} and Figure \ref{fig:edges_out}. In the employed GNN architecture, such a condition can be  conveniently enforced by integrating an additional linear constraint layer $\ell$, which was developed in our previous work \cite{chen2023learned}. Note that  $\ell$ does not have trainable parameters. It takes  %$\widetilde{d}^m_i=
 $\{\widetilde{d}^m_{ij}|j\in\mathcal{N}_{out}(i)\}$ as its input and outputs $\{d^m_{ij}|j\in\mathcal{N}_{out}(i)\}$ satisfying the conservation constraint \eqref{eq:mass}. We denote by $\bd^m = \{d^m_{ij}|e_{ij}\in \mathcal{E}\}$, and again, with a slight abuse of notations,  write $\bd^m  = \ell(\widetilde{\bd}^m)$. In summary, the decoder $Dec$ is defined as 
 \begin{equation}
 \label{eq:dec} 
 \bd^m = Dec(\bm{w}^m) := \ell\circ g(\bm{w}^m)
 \end{equation}

% After obtaining $\bd^m$, we can update the solution $\bU^{m+1}$ with \eqref{eq:sl_re}, just like the standard SL formulation, that is
% \begin{equation}\label{eq:sol-update}
%     U^{m+1}_i = \sum_{j\in N_{in}(i)}d_{ji}U^m_j,
% \end{equation}

%  we propose to integrate a constraint layer into the decoder.
% since $U_{i}^m$  contributes to computing $\{U_{j}^{m+1},j\in N_o(i)\}$ as illustrated in Figure \ref{fig:slfv_graph}. 

% \begin{figure}[htbp]
%  \centerline{\includegraphics[width=0.6\textwidth]{Concat.png}}
%   \caption{Illustration of the process of concentrating node features.}\label{fig:concat}
% \end{figure}

% The principle of the constraint layer is very simple. We only need to make the sum of the outgoing edge features of a node equal to 1, which is equivalent to \eqref{eq:conser_condition}.

After obtaining $\bd^m$, we can update the solution $\bU^{m+1}$ with  \eqref{eq:slfd}, similar as the standard SL formulation using local interpolation with a prescribed  stencil.
% \begin{equation}\label{eq:sol-update}
%     U^{m+1}_i = \sum_{j\in \mathcal{N}_{in}(i)}d_{ji}U^m_j.
% \end{equation}
Note that the stencil to update $U^{m+1}_i$ is given by $\mathcal{N}_{in}(i)$.

% Therefore, we can easily add a constraint layer after the decoder to enforce mass conservation. Figure \ref{fig:slfv_graph} presents a simple situation where $U_{i}^{m+1}$ is given by $U_{i}^{m+1} = d_{i-2,i}U^m_{i-2}+d_{i-1,i}U^m_{i-1}$. 

% ??? It should be noted that for $\forall d_i^m \in \bd^m$, it does not have a fixed shape, but depends on the specific characteristic equation and the spatial and time step size. Also, for $i\neq j$, the shape of $d_i^m$ may be different from $d_j^m$. For example, the dimensions of $d_{i-1}^m$ and $d_{i}^m$ are $2\times 1$ and $3\times 1$ respectively. If we define $N_e$ as $N_e = |\mathcal{E}|$, then $N_e = |d_1^m|+|d_2^m|+\dots+|d_N^m|$.

We remark that, within our proposed GNN-based SL framework, we favor FD discretization over the FV counterpart. Such preference stems from the complexity encountered in the high-dimensional FV setting, where determining an appropriate local stencil for updating cell averages becomes challenging due to possible severe deformation of upstream cells.

% Obviously the coefficients $\bd^m$ we get are consistent with the traditional SL scheme. 

It is worthy mentioning that for the generation of the training data, we only need to collect the coarsened solution trajectories $\bU^m$ and the normalized shifts $\bxi^m$.  Hence, we are allowed to employ any effective numerical schemes to generate data, 
%not limited to an SL FD method, 
such as the RK FD WENO method,  SL FV WENO method, among many others. The proposed end-to-end GNN-based solver can effectively learn the optimal SL FD discretization, represented by $\bd^m$, aiming to replace the most expensive component of the traditional SL methods  with a data-driven approach. 
The advantages of the proposed  solver become even more pronounced in multi-dimensional cases, since it is challenging to design a high order conservative and genuinely multi-dimensional  SL FD method.
The primary goal of our method is to simplify algorithm implementation and enhance efficiency and accuracy, while ensuring mass conservation and unconditional stability.

\begin{figure}[!htbp]
 \centerline{\includegraphics[width=0.8\textwidth]{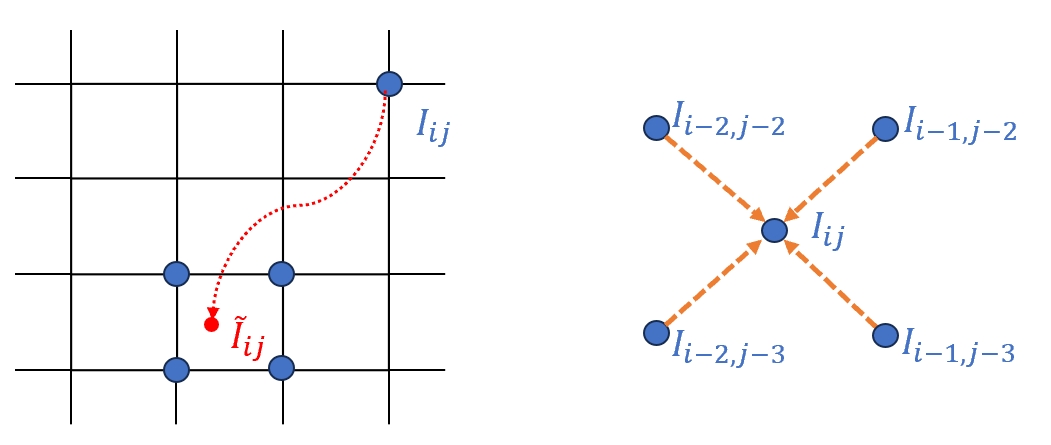}}
  \caption{Schematic illustration of the process of constructing the edges in a 2D case.}\label{fig:edges_2d}
\end{figure}

The proposed method can be generalized to the following 2D transport equation
\begin{equation}\label{eq:2d}
u_t + \nabla\cdot (\bv(x,y,t) u) = 0,\quad (x,y)\in\Omega, 
\end{equation}
where $\bv$ denotes the velocity field $\bv(x,y,t) = (a(x,y,t),b(x,y,t))$. The associated characteristic system is given by
\begin{equation}
\label{eq:characteristic2d}
  \begin{cases}
  \frac{dx(t)}{dt} &= a(x(t),y(t),t),\\
   \frac{dy(t)}{dt} &= b(x(t),y(t),t).\\
  \end{cases}
\end{equation}
The rectangular domain $\Omega$ is discretized using a set of uniformly spaced tensor product grid points $\bigcup_{ij} I_{ij}$. Denote by $U_{ij}^m$ the approximate point value of the solution $u$ at the grid point $I_{ij}$ and time step $t^m$. Similar as  the 1D case \eqref{eq:shift},  by solving the characteristic system \eqref{eq:characteristic2d}, we define the normalized shifts $\xi^m_{ij}$ and  $\eta^m_{ij}$ at each grid point $I_{ij}$ and time step $t^m$ in $x$- and $y$- directions, respectively. 

In the 2D case, we still employ the Encode-Process-Decode framework, with slight modifications in  the encoder. Specifically, the encoder $Enc$ is defined as 
\begin{equation}
\bm{h}^{m,(0)} = Enc(\bU^m,\bxi^m,\bm{\eta}^m),
\end{equation}
where the inputs $\bU^m,\bm{\xi}^m$, and $\bm{\eta}^m$ denote the collections  $\{U_{ij}^m\}$, $\{\xi_{ij}^m\}$, and $\{\eta_{ij}^m\}$, respectively. $Enc$ takes 3-channel node features as input and computes node embedding vectors. In particular, $Enc$ is constructed by staking a sequence of 2D convolutional layers together with nonlinear activation functions. Meanwhile, the structure of the GNN used in the processor is exactly the same as that in the 1D case. In particular, we  model the grid as a graph $\mathcal{G}=(\mathcal{V},\mathcal{E})$, where the node set $ \mathcal{V}$ comprises all the grid points ${I_{ij}}$, and $e_{i_1j_1,i_2j_2} \in \mathcal{E}$ represents the directed edge that connects the node $I_{i_1j_1}$ to the node $I_{i_2j_2}$. To construct the edge set $\mathcal{E}$, we first locate the upstream point $\widetilde{I}_{ij}$ for each node $I_{ij} \in \mathcal{V}$. Next, we select a local stencil centered at the upstream point $\widetilde{I}_{ij}$, such as the four-point stencil composed of the four nearest grid points $I_{i_1j_1}$, $I_{i_2j_2}$, $I_{i_3j_3}$, and $I_{i_4j_4}$. We then create four directed edges from these nodes to the node $I_{ij}$, as illustrated in Figure \ref{fig:edges_2d}. Lastly, the decoder is identical to that is used in the 1D case, including an MLP and a constraint layer to enforce mass conservation. Once the coefficients $\bd^m$ are determined, the solution  is updated to the next time step as in the 1D case \eqref{eq:slfd}. 
\\

\noindent{\bf Computational Complexity.} We provide the forward computational complexity of each block in the proposed Encode-Process-Decode framework. $|\mathcal{V}|$ denotes the number of grid points or nodes, and  $|\mathcal{E}|$ denotes the number of edges in the graph. For the CNN encoder, the computational complexity is $O(|\mathcal{V}| c_K c_{in} c_{out} )$ for each convolutional layer, where $c_K,\, c_{in},$ and $c_{out}$ are the kernel size, the number of input channels, and the number of output channels, respectively.  In the 2D case, the computational complexity becomes $O(|\mathcal{V}| c_K^2 c_{in}  c_{out} )$ with the kernel size being $c_K \times c_K$. For the GNN processor, the computational complexity of one graph attention layer is $O(|\mathcal{V}|FF'+|\mathcal{E}|F' )$, where $F$ and $F'$ denote the number of input features and the number of output features of the graph node, respectively. For the two components of the decoder, the complexity of each linear layer in the MLP is $O(F'|\mathcal{E}|)$, and the complexity of the constraint layer is $O(|\mathcal{E}|)$. Note that  the size of the stencil $\mathcal{N}_{in}(i)$ is fixed for each node $I_i$, implying that $|\mathcal{E}| = O(|\mathcal{V}|)$.

\subsection{The Vlasov-Poisson system}\label{sec:vp}
In this subsection, we extend the proposed algorithm to the nonlinear VP system,  addressing the challenges posed by  its inherent nonlinearity.

The VP system is a fundamental model in plasma physics that describes interactions between charged particles through self-consistent electrostatic fields, modeled by Poisson’s equation. The dimensionless VP equation under 1D in space and 1D in velocity (1D1V) setting is given by
\begin{align}
    &f_t + vf_x + E(x,t)f_v = 0,\quad (x,v)\in\Omega_x\times\Omega_v,\label{eq:vlasov}\\
    &E_x = \rho - 1,\quad \rho(x,t) = \int_{\Omega_v} f(x,v,t) dv,\label{eq:poisson}
\end{align}
where $f(x,v,t)$ is the probability distribution function of electrons at position $x$ with velocity $v$ at time $t$, and  $\rho$ denotes the macroscopic density. $\Omega_x$ denotes the physical domain, while  $\Omega_v$ represents the velocity domain. In addition, we assume a uniform background of fixed ions and consider  periodic boundary conditions in the $x$-direction and  zero boundary conditions in the $v$-direction. It is worth mentioning that in practice $\Omega_v$ is truncated to be finite and is taken large enough so that the solution $f\approx 0$ at $\partial \Omega_v$.

 % under a self-consistent electrostatic field $E$.
%We assume that the ions are fixed and form a neutralizing background. The VP system is known as a reduced model of the full Vlasov-Maxwell system in the zero-magnetic limit and plays a vital role in plasma physics.

Unlike \eqref{eq:transport1d}, the Vlasov equation \eqref{eq:vlasov} is a nonlinear transport equation. The nonlinearity introduces additional challenges in accurately approximating its characteristic equations:
\begin{equation}\label{eq:vpchar}
\begin{cases}
 \displaystyle  \frac{d x(t)}{d t} &= v(t),\\[3mm]
  \displaystyle  \frac{d v(t)}{d t} &= E(x(t),t).
\end{cases}
\end{equation}
% Consequently, designing an SL scheme becomes more complex for the VP system.
To circumvent the difficulty, the splitting approach was introduced in \cite{cheng1976integration,sonnendrucker1999semi}, which decouples the nonlinear Vlasov equation into several linear equations in lower dimensions.
 However, the inherent splitting error may greatly compromise the overall accuracy of the approximation.
Recently, a non-splitting SL methodology has been introduced in \cite{celledoni2003commutator} for the VP system, utilizing  the commutator-free RKEI. By employing the RKEI method, the VP system is decomposed into a series of linearized transport equations with frozen coefficients \cite{zheng2022fourth,cai2021high}, which can be effectively solved using the proposed GNN-based SL FD scheme. With RKEI, our GNN-based SL FD VP solver eliminates splitting errors while maintaining all the advantages in the linear case.

To introduce the algorithm, we begin by considering a uniform partition of the domain $\Omega_x\times\Omega_v$ with $N_x\times N_v$ grid points, i.e., $\Omega_x\times\Omega_v=\bigcup_{ij}I_{ij}$, where the grid point $I_{ij}$ has coordinate $(x_i,v_j)$. The mesh sizes in the $x$- and $v$- directions are represented by $h_x$ and $h_v$, respectively. Let  $f^m_{ij}$ denote the numerical approximation of $f$ at the grid point $I_{ij}$ at time level $t^m$ and $E_{i}^m$ be the numerical solution of $E$ at the grid point $I_i$ in the $x$- direction. The collections $\{E^m_{i}\}$ and $\{f^m_{ij}\}$ are denoted by $\bE^m$ and $\bF^m$, respectively.
%  
%the cell average of the density $\rho$, denoted by $\rho^m_i$, can be computed by $$\rho^m_i =h_v\sum_{j} f_{ij}^m,$$ and 
The electric field $\bE^m$ is solved from Poisson's equation \eqref{eq:poisson}.
\begin{table}[!htbp]
\begin{minipage}{.5\linewidth}
    \centering

    \caption{First-order RKEI}
    \label{tab:1strkei}

    \medskip

\begin{tabular}{l|l}
0 & 0 \\
\hline & 1
\end{tabular}
\end{minipage}\hfill
\begin{minipage}{.5\linewidth}
    \centering

    \caption{Second-order RKEI}
    \label{tab:2ndrkei}

    \medskip

\begin{tabular}{l|ll}
0 & & \\
$\frac{1}{2}$ & $\frac{1}{2}$ & 0 \\
\hline & 0 & 1
\end{tabular}
\end{minipage}
\end{table}

Our approach employs the RKEI technique to circumvent the difficulties associated with accurately tracking the characteristics for the VP system. The RKEI method is represented by a Butcher tableau, which specifies the coefficients for the integration steps. The accuracy of the method is determined by order conditions. The simplest RKEI method is given by the Butcher tableau \ref{tab:1strkei}, which is first order accurate \cite{celledoni2003commutator}. With the first-order RKEI, we can develop a GNN-based SL FD scheme for the VP system as follows: 
\begin{itemize}
    \item[(1)] Compute $\bE^m$ from Poisson equation. 
    \item[(2)] Linearize the Vlasov equation with the fixed electric field  $\bE^m$ and solve the characteristic equations \eqref{eq:vpchar}  to obtain the normalized shifts.
    \item[(3)] Use the GNN-based SL FD method, together with $\bF^m$ and the normalized shifts, to predict $\bF^{m+1}$, as with the linear case.
\end{itemize}
The procedure  can be summarized as
\begin{equation}
\label{eq:1stvp}
\bF^{m+1} = GNNSL(\bE^m,\Delta t)\bF^{m},
\end{equation}
where $GNNSL(\bE^m,\Delta t)$ denotes the proposed GNN-based SL FD evolution operator for the Vlasov equation with the fixed electric field  $\bE^m$ and time step $\Delta t$.
However, the low order temporal accuracy of the first order RKEI  can severely limit its  generalization capability, as observed in \cite{chen2023multi}. 
Meanwhile, we may enhance the performance by employing a second order RKEI \cite{celledoni2003commutator}, represented by the Butcher tableau \ref{tab:2ndrkei}.
The corresponding GNN-based SL FD algorithm for the VP system can be summarized as 
\begin{equation}
\label{eq:2ndvp}
\begin{cases}
        \bF^{m,*} = GNNSL(\bE^m,\frac12\Delta t)\bF^{m},  \\
    \bF^{m+1} = GNNSL(\bE^{m,*},\Delta t)\bF^{m}, 
\end{cases} 
\end{equation}
where $\bE^m$ and $\bE^{m,*}$ are determined from $\bF^{m}$ and $\bF^{m,*}$, respectively. The second order scheme involves an intermediate stage $\bF^{m,*}$ and requires two applications of the proposed GNN-based SL FD evolution operator. Note that in the process of $\bF^{m} \to \bF^{m+1}$, the dynamical graph is updated twice. In particular, the edge set is first determined by $\bE^m$, and then updated by $\bE^{m,*}$. As with the linear case,  we are allowed to take extra large time step evolution for simulating the nonlinear VP system due to the proposed dynamical GNN architecture.
Furthermore, it is numerically demonstrated that the second order method \eqref{eq:2ndvp} exhibits improved  generalization capabilities and achieves higher accuracy compared to the first counterpart \eqref{eq:1stvp}. Hence, we only present the numerical results obtained by the second order RKEI in the next section.   
% Moreover, all the desired properties of the method for solving linear equations, such as mass conservation, are preserved when applied to the VP system. 
Such a neural VP solver achieves a level of accuracy that surpasses traditional numerical algorithms, such as the popular SL FV WENO scheme and the Eulerian RK WENO scheme, with comparable mesh resolution. 

It is worth mentioning that one-step training, as described above, may suffer slightly weaker generalization. One the other hand, unrolling the training process over multiple time steps can improve the accuracy and stability at the cost of increased training difficulty, as discussed in \cite{brandstetter2021message}. One-step training is sufficient for linear transport equations, while the training is unrolled with eight time steps for the nonlinear VP system.
% \subsection{Discussion}

% We end this section with some comments on the proposed method in the following:
% \begin{itemize}
%     \item xxx
% \end{itemize}

\section{Numerical results}
\label{sec:experiments}
In this section, we carry out a series of numerical examples to demonstrate the performance of the proposed GNN-based SL FD scheme for simulating various benchmark 1D and 2D linear transport equations together with the nonlinear 1D1V VP system. Noteworthy, the performance of the proposed method depends on the structure of each block within the Encoder-Processor-Decoder framework, as illustrated in Figure \ref{fig:network}. It includes the choice of hyperparameters. For simplicity, the numerical results are presented using  default settings. For the linear transport equations, the CNN encoder is configured with six convolutional layers, each equipped with 32 filters, and a kernel size of five for both the 1D and 2D cases.  
For the GNN processor, we utilize two graph attention layers, each with an output feature dimension of 32 and four attention heads. Lastly, the MLP used in the decoder features  one hidden layer with 256 neurons.  For the nonlinear VP system, the CNN encoder is configured with nine convolutional layers, each containing 32 filters, and  a kernel size of five. The configurations of the processor and decoder are identical to those used in the linear case.  We employ ELU \cite{clevert2016fast} as the activation function and Adam \cite{Ilya_fix_2017} as the optimizer in the implementation.

As discussed in Section \ref{sec:main}, we can employ any accurate and reliable numerical scheme to generate the training data. In this paper, for the linear transport equations, we use the fifth-order WENO (WENO5)  FD method \cite{jiang1996efficient}, combined with the third order strong-stability-preserving RK (SSPRK3) time integrator \cite{gottlieb2001strong}. For the nonlinear VP system, we employ the conservative fifth-order SL FV WENO scheme coupled with a fourth-order RKEI \cite{zheng2022fourth}. The training data is generated by coarsening high-resolution solution trajectories onto a low resolution mesh. For all the test examples, we mainly report the results by the proposed GNN-based SL FD method and the WENO5 method with the same mesh resolution, together with the reference solutions for comparison.  In all the plots reported below, ``GNN" denotes the proposed method. 
% \lipsum[50]

\subsection{One-dimensional transport equations}
In this subsection, we present numerical results for simulating 1D transport equations. 
%In this case, the node feature for node $i$ is $x_i = (U_i^m,\xi_i^m)\in \mathbb{R}^2$.

%The spatial domain is taken as $[0,1]$. We consider a one-dimensional grid with 32 grid points.
\begin{example}\label{square}
In this example, we consider the following advection equation with a constant coefficient
\begin{equation}\label{eq:const1d}
u_t + u_x = 0,\quad x\in [0,1],
\end{equation}
and periodic boundary conditions are imposed. %we consider the square waves.
% The initial condition is taken as square waves and the solution profile is shifted by a constant distance per time step, determined by the CFL number.
\end{example}

The training data is generated by coarsening 30 high-resolution solution trajectories on a 256-cell grid by a factor of 8. The initial condition for each trajectory is a square wave with height randomly sampled from $[0.1,1]$ and width from $[0.2,0.4]$. Each trajectory contains 20 sequential time steps, and the CFL numbers are chosen within the range of $[6,10.2]$. To test, we randomly generate square functions with heights and widths within the same range as initial conditions. For comparison, the reference solution is produced using WENO5 method on a high-resolution 256-cell grid and then down-sampled to a coarse grid of 32 cells, the procedure of which is the same as generating the training data.

\begin{figure}[!htbp]
 \centering   
\includegraphics[width=1.00\textwidth]{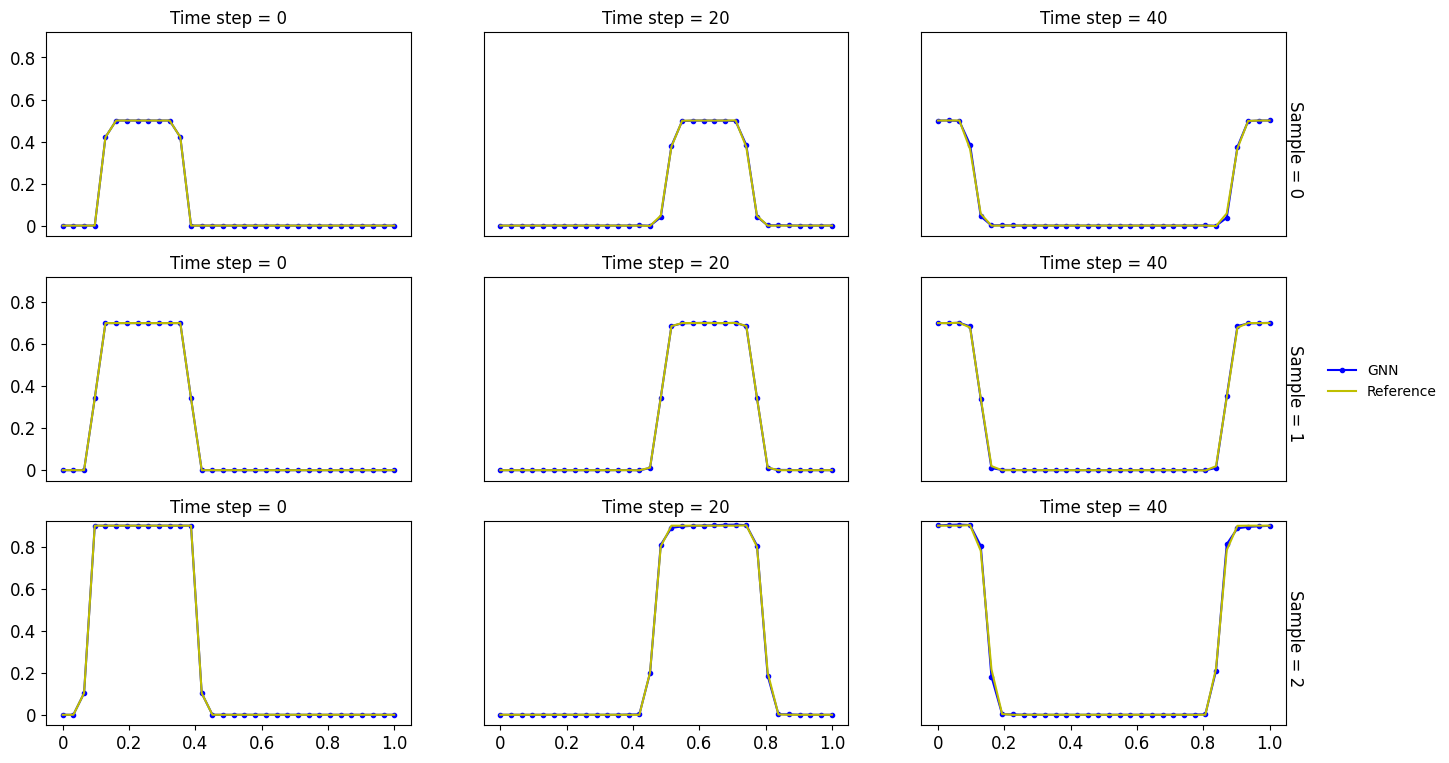}
\caption{Three test samples for square waves in Example \ref{square}.  CFL=10.2.}\label{fig:square_three}
\end{figure}

\begin{figure}[!htbp]
 \centering
     \subfigure[]{\includegraphics[width=0.37\textwidth,height=0.36\textwidth]{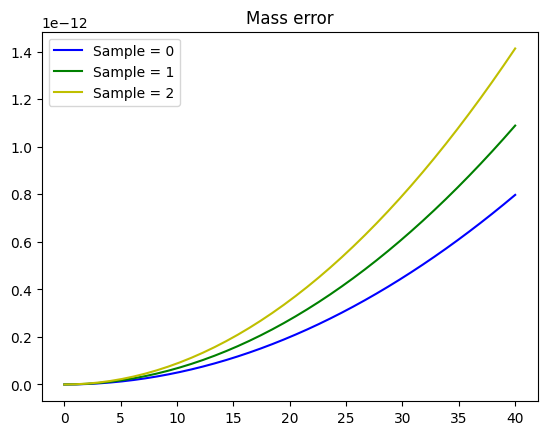}}
    \subfigure[]{\includegraphics[width=0.39\textwidth,height=0.355\textwidth]{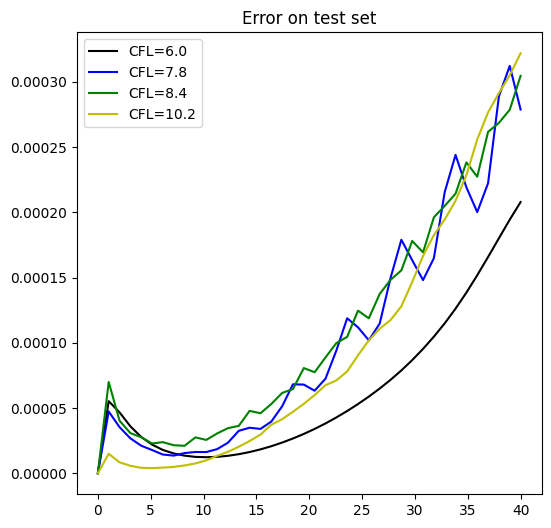}}
	\caption{Example \ref{square}. Time histories of the  deviation of total mass for three test samples with CFL = 10.2 (a), and  time histories of errors with different CFLs (b).}\label{fig:square_mass}
\end{figure}

Figure \ref{fig:square_three} displays plots of three test samples during forward integration at different time instances, using a CFL number of 10.2, which significantly exceeds the CFL limit for the  WENO5 scheme as well as our previous ML-based SL FV method \cite{chen2023learned}. It is observed that the proposed GNN-based approach demonstrates superior shock resolution, with sharp shock transitions and no spurious oscillations, and the results stay close to the reference solution. In Figure \ref{fig:square_mass} (a), we presents the time evolution of the total mass deviation of the three test solution trajectories generated by our method. Clearly, the total mass is conserved to the machine precision. 

\begin{figure}[!htbp]
    \centerline{\includegraphics[width=0.7\textwidth]{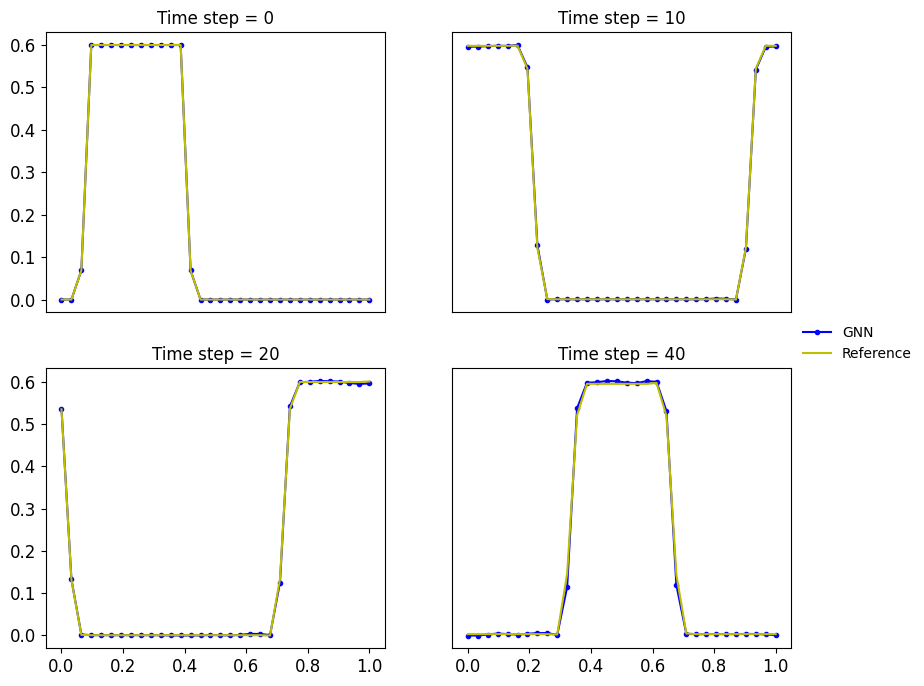}}
	\caption{One test example for square waves in Example \ref{square}. CFL=9.0. The CFL number is not contained in the training data. }\label{fig:square_cfl_test}
\end{figure}

% Since the training data contains solution trajectories for several CFL numbers, 
We investigate the  performance of the proposed method  with different CFL numbers used in the training data, i.e., different time step sizes, and display the result in Figure \ref{fig:square_mass} (b). We can observe that the errors by our method with different CFLs are almost of the same magnitude over time. Furthermore, in Figure \ref{fig:square_cfl_test}, we also present one test example obtained from forward integration at several instances of time with CFL=9, which is not used in the training data. The performance is comparable to that observed in Figure \ref{fig:square_three}. We remark that the proposed method permits the use of arbitrary CFL numbers within the range used to generate the training data. When the CFL number exceeds this range, there is a deterioration in performance.

% \begin{figure}[!htbp]
%     \centerline{\includegraphics[width=0.5\textwidth]{square_cfl.png}}
% 	\caption{Errors by the proposed method with
% different CFLs.}\label{fig:square_cfl}
% \end{figure} 

\begin{example}\label{triangle}
In this example, we consider the advection equation \eqref{eq:const1d} with a more complicated solution profile consisting of triangle and square waves. 
\end{example}

Similar to the previous example, we generate 30 high-resolution solution trajectories, each consisting of one triangle and one square waves with heights randomly sampled from $[0.2,0.8]$ and widths from $[0.2,0.3]$ over the 256-cell grid. By coarsening these solution trajectories by a factor of 8, we obtain our training data. Each solution trajectory in the training data set contains 20 time steps with the CFL numbers ranging in $[6,10.2]$. For testing, we generate initial conditions for test data within the same ranges of width and height. To evaluate the performance of the proposed method, we calculate the ground-truth reference solution using WENO5 on a high-resolution grid with 256 cells and then reduce the resolution by a factor of 8 to a coarse grid of 32 cells.

\begin{figure}[!htbp]
    \centerline{\includegraphics[width=0.95\textwidth]{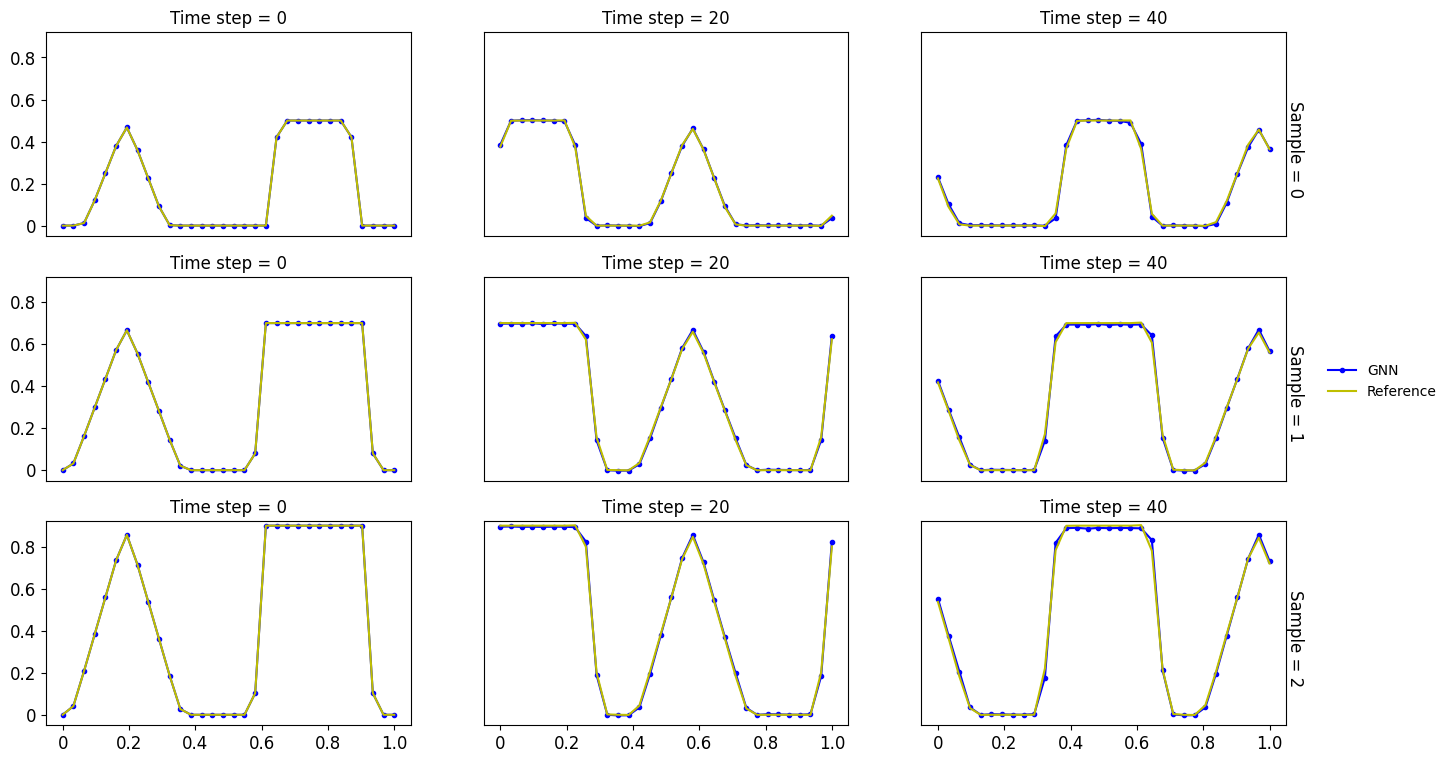}}
	\caption{Numerical solutions of three test samples for advection of triangle and square waves in Example \ref{triangle}.  CFL=10.2.}\label{fig:tri_three}
\end{figure}

\begin{figure}[!htbp]
 \centering
     \subfigure[]{\includegraphics[width=0.37\textwidth,height=0.36\textwidth]{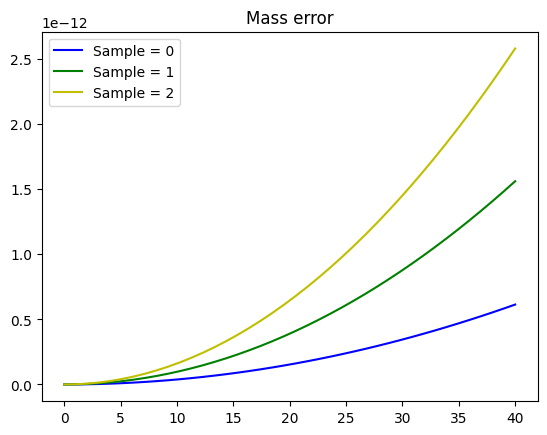}}
     \subfigure[]{\includegraphics[width=0.39\textwidth,height=0.36\textwidth]{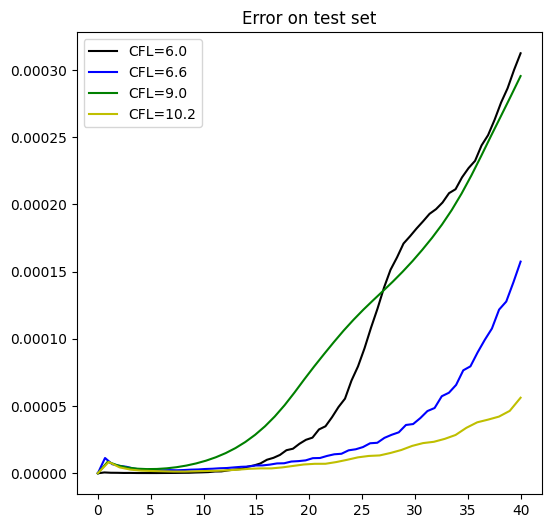}}
	\caption{Example \ref{triangle}. Time histories of the  deviation of total mass for three test samples with CFL = 10.2 (a), and  time histories of errors with different CFLs (b).}\label{fig:tri_mass}
\end{figure}

In Figure \ref{fig:tri_three}, we present three test examples at several instances of time during forward integration with CFL=10.2. Similar to the previous example, the proposed GNN-based solver achieves superior shock resolution, and the simulation results stay close to the reference solution over time. Our method is also mass conservative up to machine precision as demonstrated in Figure \ref{fig:tri_mass} (a). We further validate our solver with different CFL numbers, and plot the time histories of errors in Figure \ref{fig:tri_mass} (b). It is observed that the method using a larger CFL number results in a smaller error.
% In addition, we test the solver's ability to generalize with arbitrary CFL numbers. Figure \ref{fig:tri_cfl_test} plots one test example during forward integration at several time instances with CFL=7.8, which is not included in the training data.

% \begin{figure}[!htbp]
%     \centerline{\includegraphics[width=0.8\textwidth]{tri_cfl.png}}
% 	\caption{Errors by the proposed method with
% different CFLs.}\label{fig:tri_cfl}
% \end{figure} 

% \begin{figure}[!htbp]
%     \centerline{\includegraphics[width=0.8\textwidth]{tri_cfl_test.png}}
% 	\caption{One test example for square waves in Example \ref{triangle}. CFL=7.8. The CFL number is not contained in the training data. }\label{fig:tri_cfl_test}
% \end{figure} 

\begin{example}
\label{eg:variable}
In this example, we simulate the following 1D advection equation with a variable coefficient
\begin{equation}\label{eq:variable-velocity}
    u_t+(\sin(x+t)u)_x=0,\quad x\in[0,2\pi],
\end{equation}
subject to periodic boundary conditions. This example is more challenging than the previous two examples. In addition to  pure shift, solution profiles also deform  gradually over time and exhibit more complex structures.

For this example,  as time $t$ increases, the solution will gradually concentrate mass at a single point and converge towards a $\delta$ function. In particular, when $t$ is approximately greater than 3.5,  the mesh cannot provide adequate resolution for such a singular structure, resulting in irrelevant results. In addition, if we choose  CFL=9, then  the time step size $\Delta t$ is $\frac{9\pi}{16}\approx 1.77$. Given such large time step utilized, we have to limit the simulation to only one or two steps. 
\end{example}

To generate the training data, we coarsen 90 solution trajectories on a 256-cell grid by a factor of 8 with each trajectory consisting of 2 sequential time steps. The CFL numbers are chosen within the range of $[5,9]$. As in the first example, the initial condition is a step function with heights randomly sampled from the range $[0.1,1]$ and widths sampled from the range $[2.5,3.5]$. In addition, the center of each square function is randomly sampled from the whole domain $[0,2\pi]$.  Again, the reference solution is generated by WENO5 over the 256-cell grid  and down-sampled to the coarse grid of 32 cells.  

\begin{figure}[!htbp]
    \centerline{\includegraphics[width=0.95\textwidth]{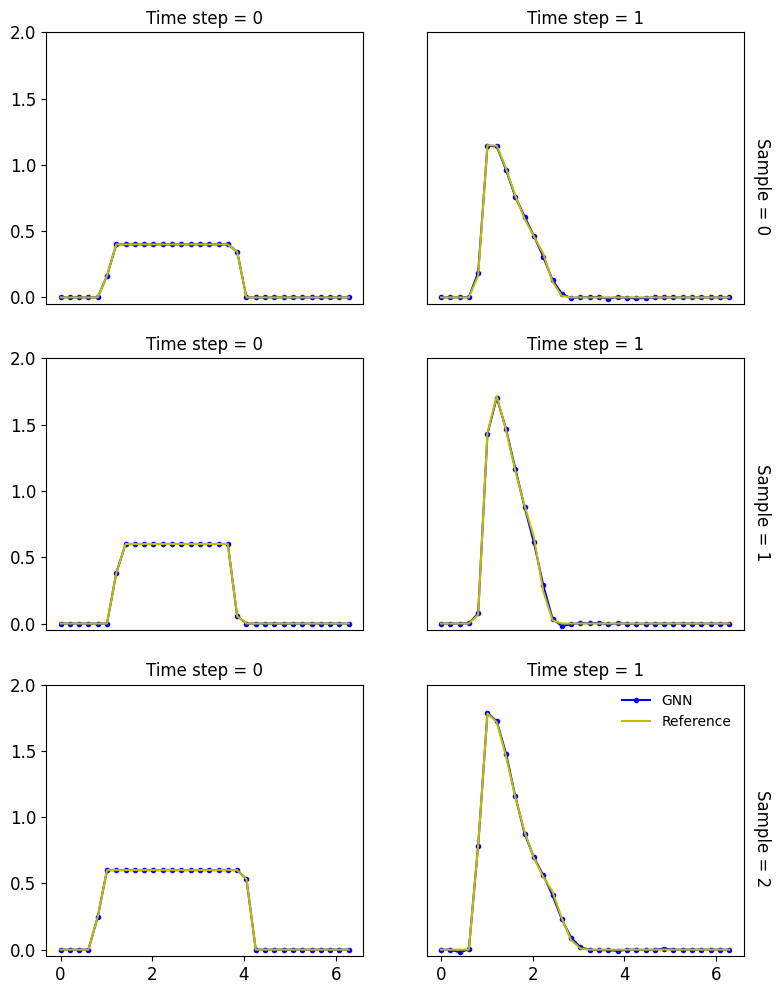}}
	\caption{Numerical solutions of three test samples for the transport equation with a variable coefficient in Example \ref{eg:variable}.  CFL=9.0.}\label{fig:var_three}
\end{figure}

\begin{figure}[!htbp]
 \centering
     \subfigure[]{\includegraphics[width=0.37\textwidth,height=0.36\textwidth]{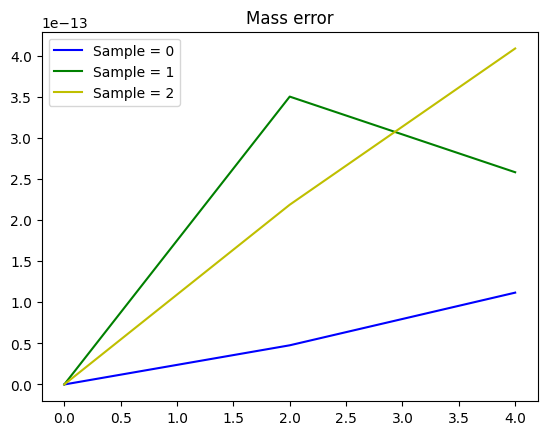}}
     \subfigure[]{\includegraphics[width=0.39\textwidth,height=0.36\textwidth]{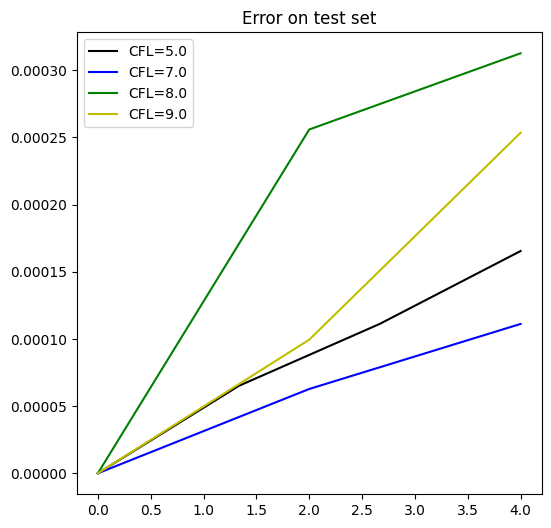}}
	\caption{Example \ref{eg:variable}. Time histories of the  deviation of total mass for three test samples with CFL = 9.0 (a), and  time histories of errors with different CFLs (b).}\label{fig:var_mass}
\end{figure}

\begin{figure}[!htbp]
    \centerline{\includegraphics[width=0.95\textwidth]{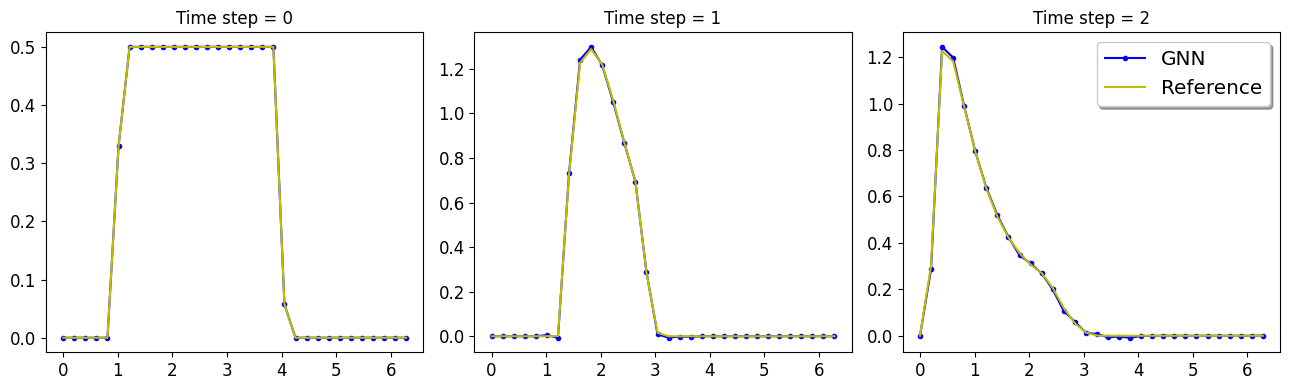}}
	\caption{One test example for the transport equation with a variable coefficient in Example \ref{eg:variable}. CFL=6.0. The CFL number is not contained in the training data. }\label{fig:var_cfl_test}
\end{figure}

Figure \ref{fig:var_three} shows three test samples, each updated in a single step with CFL=9. It is observed that our solver can accurately resolve singular solution structures with sharp transitions, despite the use of a very large time step size $\Delta t \approx1.77$. As demonstrated in Figure \ref{fig:var_mass} (a), the proposed method can conserve the total mass up to machine precision. We further report the time histories of errors for the method with different CFL numbers in Figure \ref{fig:var_mass} (b). We can observe that the errors are comparable across different CFLs. Moreover, Figure \ref{fig:var_cfl_test} presents one test example during forward integration in two steps with CFL=6, which is not used during training, yet high quality numerical results are observed.
This indicates that our method can utilize arbitrary CFLs within the range used to generate the training data.

% \begin{figure}[!htbp]
%  \centerline{\includegraphics[width=0.45\textwidth]{var_mass.png}}
% 	\caption{Time histories of the  deviation of total mass of three test samples for the transport equation with a variable coefficient in Example \ref{eg:variable}. CFL=10.2.}\label{fig:var_mass}
% \end{figure} 

% \begin{figure}[!htbp]
%     \centerline{\includegraphics[width=0.8\textwidth]{var_cfl.png}}
% 	\caption{Errors by the proposed method with
% different CFLs.}\label{fig:var_cfl}
% \end{figure} 

\subsection{Two-dimensional transport equations}
In this subsection, we present the numerical results for simulating several 2D benchmark advection problems. %In the 2d case, the node feature for node $I_{ij}$ is $x_{ij}= (U_{ij}^m,\xi_{ij}^m,\eta_{ij}^m)\in \mathbb{R}^3$.

\begin{example}
\label{eg:2dlinear}
 In this example, we consider the following constant-coefficient 2D transport equation 
$$
u_t + u_x + u_y= 0,\quad (x,y)\in[-1,1]^2,
$$
with periodic boundary conditions.
\end{example}

We generate the training data by coarsening 30 high-resolution solution trajectories over a
$256\times 256$-cell grid by a factor of 8 in each dimension. Each trajectory contains 15 sequential time steps with CFL=10.2. The initial condition for each trajectory is a square wave with height randomly
sampled from $[0.5,1]$ and width from $[0.3,0.5]$. After training, we test the performance of the solver for the problem with initial conditions sampled from the same ranges of width and height. The reference solution is generated by WENO5 over the $256\times 256$ cells and down-sampled to the coarsen grid of $32\times 32$ cells.

\begin{figure}[!htbp]
    \centerline{\includegraphics[width=0.9\textwidth]{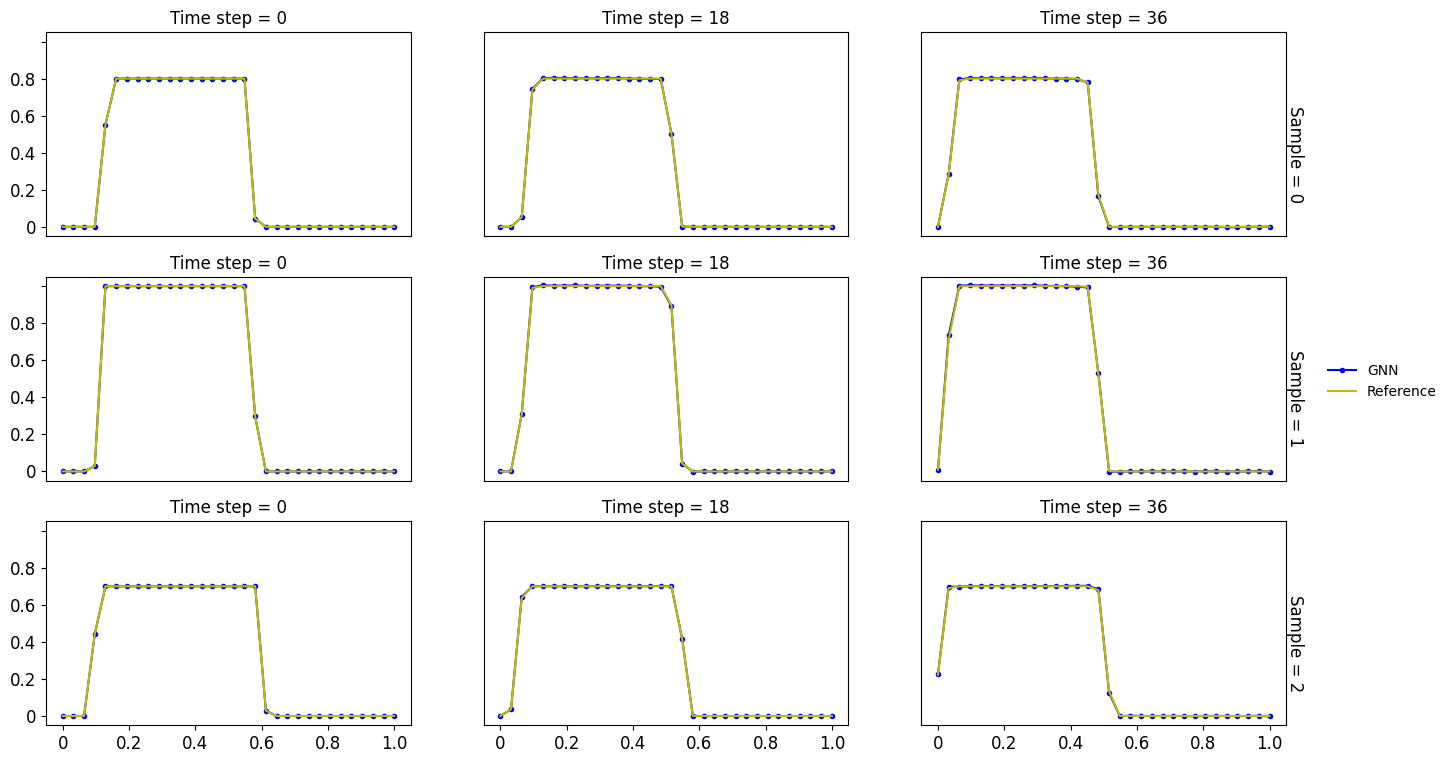}}
	\caption{1D cuts at $y=x$ of three test samples for 2D transport equation with constant coefficients in Example \ref{eg:2dlinear}.  CFL=10.2.}\label{fig:2d_three_curve}
\end{figure}

\begin{figure}[!htbp]
    \centerline{\includegraphics[width=0.9\textwidth]{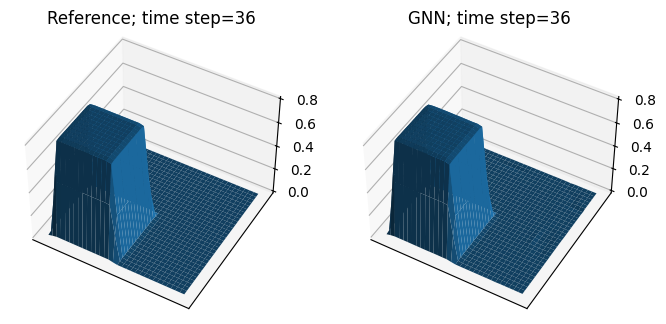}}
	\caption{One test samples for the 2D transport equation with constant coefficients in Example \ref{eg:2dlinear}. 2D plots of the solutions at time step 36. CFL=10.2.}\label{fig:2d_curve}
\end{figure}

\begin{figure}[!htbp]
 \centerline{\includegraphics[width=0.45\textwidth]{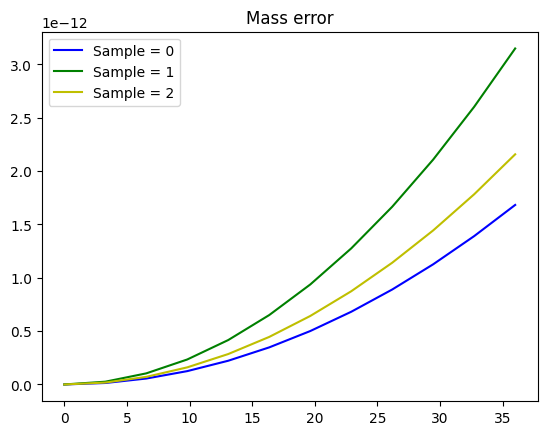}}
	\caption{Time histories of the  deviation of total mass of three test samples for the 2D transport equation with a constant coefficient in Example \ref{eg:2dlinear}. CFL=10.2.}\label{fig:2d_mass}
\end{figure}

In Figure \ref{fig:2d_three_curve}, we present 1D cuts of solutions at $y=x$ for three test examples plotted at different time instances during forward integration with CFL=10.2. For a more effective comparison, we also provide the 2D plots of the solutions at time step 36 in Figure \ref{fig:2d_curve}. It can be observed that our solver can resolve discontinuities sharply without introducing spurious oscillations, despite the use of such a large CFL number. Figure \ref{fig:2d_mass} shows the time evolution of the total mass deviation of three solution trajectories generated by our solver. Evidently, the total mass is conservative up to the machine precision.


\begin{example}\label{defor}
In this example, we simulate a 2D linear  deformation flow problem proposed in \cite{leveque_high-resolution_1996}, governed by the following  transport equation
  \begin{equation}\label{eq:2d:defor}
u_t + (a(x,y,t)u)_x+(b(x,y,t)u)_y = 0,\quad (x,y)\in[0,1]^2. 
\end{equation}
The velocity field is a periodic swirling flow \begin{equation}\label{eq:deformational flow}
\begin{aligned}
    a(x,y,t) &= \sin^2(\pi x)\sin(2\pi y)\cos(\pi t/T),\\
    b(x,y,t) &= -\sin^2(\pi y)\sin(2\pi x)\cos(\pi t/T),
\end{aligned}
\end{equation}
where $T$ is a constant.
It is a widely recognized benchmark test for numerical transport solvers. It exhibits distinct dynamic properties in which the solution profile deforms over time in response to the flow. The direction of the flow reverses att $t=T/2$, and the solution returns to its initial state at $t=T$, completing a full cycle of the evolution.
\end{example}

We set $T=2$ and choose the initial condition to be a cosine bell centered at $[c_x,c_y]$
\begin{equation}\label{eq:deformation-initial}
\begin{aligned}
    u(x,y) &= \frac{1}{2}[1+\cos (\pi r)],\\
    r(x,y) &= \min\left[1,r_0\sqrt{(x-c_x)^2+(y-c_y)^2}\right],
\end{aligned} 
\end{equation} 
where $r_0$ determines the radius of the cosine bell. To generate the training data, we initialize 90 trajectories with $c_x,c_y$ randomly sampled from $[0.25,0.75]$ and $r_0$ randomly sampled from $[4,6]$ using a high-resolution mesh of $256\times256$ cells. Then we coarsen these trajectories by a factor of 8 in each dimension. Each solution trajectory contains a sequence of time steps from $t=0$ to $t=T$, with CFL=10.2. To test, we choose initial conditions  sampled from the same distribution as the training data. For the purpose of comparison, the reference solution is generated with the same approach used to generate the training data.

\begin{figure}[!htbp]
    \centerline{\includegraphics[width=0.9\textwidth]{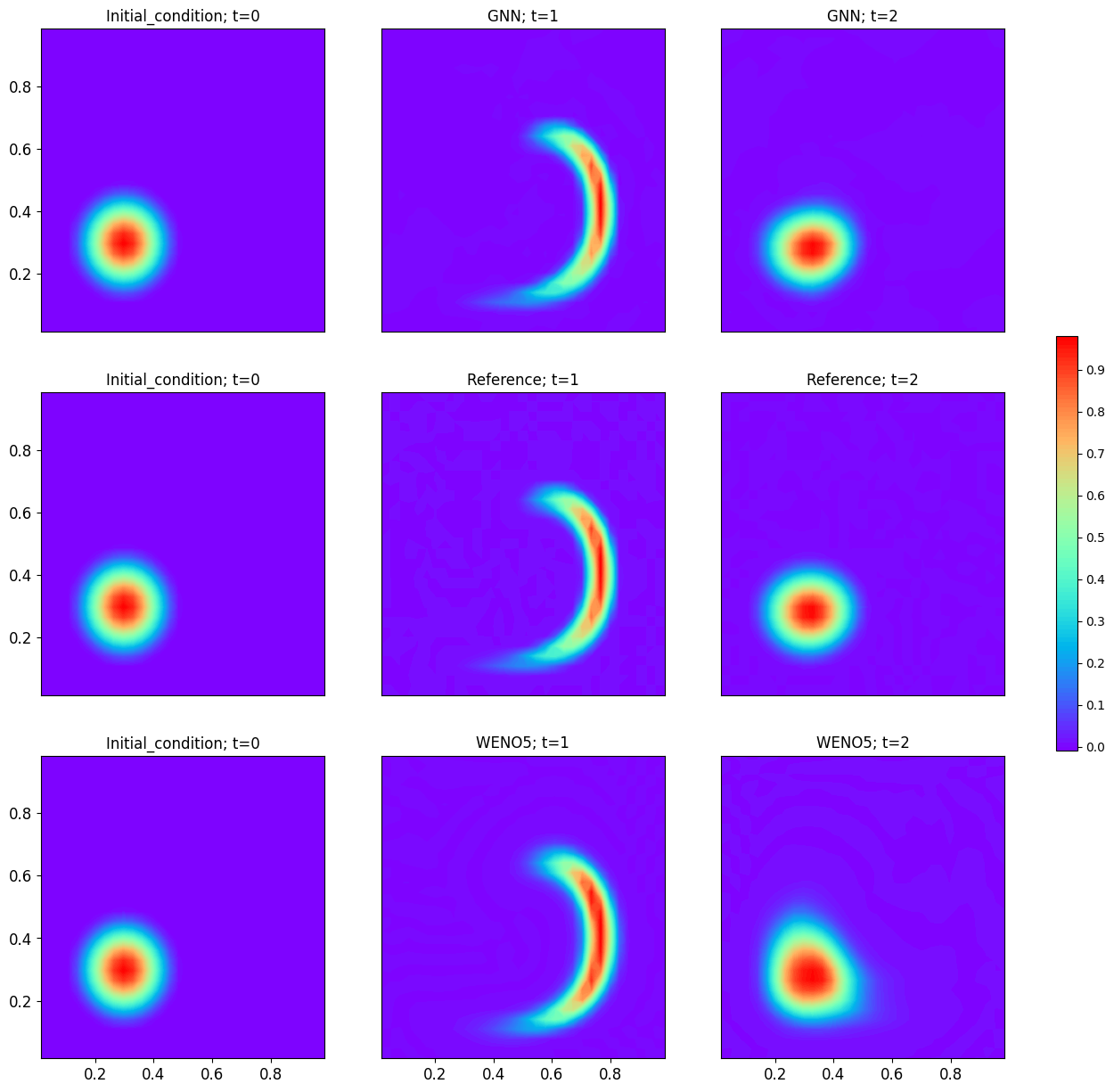}}
	\caption{Contour plots of the numerical solutions of the 2D deformational flow at $t=0,\,1,\,2$  in Example \ref{defor} for one test sample. CFL=10.2 for our method, CFL=0.6 for WENO5. $r_0=5,c_x=0.3,c_y=0.3$. This configuration is not contained in the training data.}\label{fig:defor_test}
\end{figure}

\begin{figure}[!htbp]
 \centerline{\includegraphics[width=0.45\textwidth]{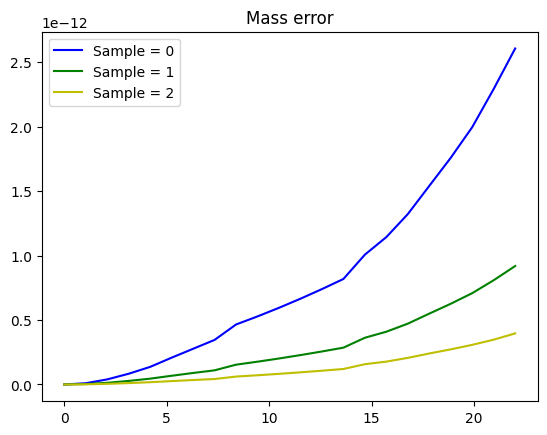}}
	\caption{Time evolution of the  deviation of total mass for three test samples in Example \ref{eg:2dlinear}. CFL=10.2.}\label{fig:deform_mass}
\end{figure}

Figure \ref{fig:defor_test} shows the contour plots of the numerical solutions computed by our methods and WENO5 together with the reference solution for one test example with $r_0=5,c_x=0.3,c_y=0.3$, a configuration  that is not included in the training set. The CFL numbers are chosen as 10.2 and 0.6 for our method and WENO5, respectively. It is observed that the solution is significantly deformed at $t=T/2$ and returns to its initial state at $t=T$. Our solver can accurately capture deformations and restore the initial profile. However, the solution obtained using WENO5 noticeably deviates from the reference solution due to a large amount of numerical diffusion. As demonstrated in Figure \ref{fig:deform_mass}, our method can also conserve the total mass up to machine precision.

\begin{figure}[!htbp]
    \centerline{\includegraphics[width=0.9\textwidth]{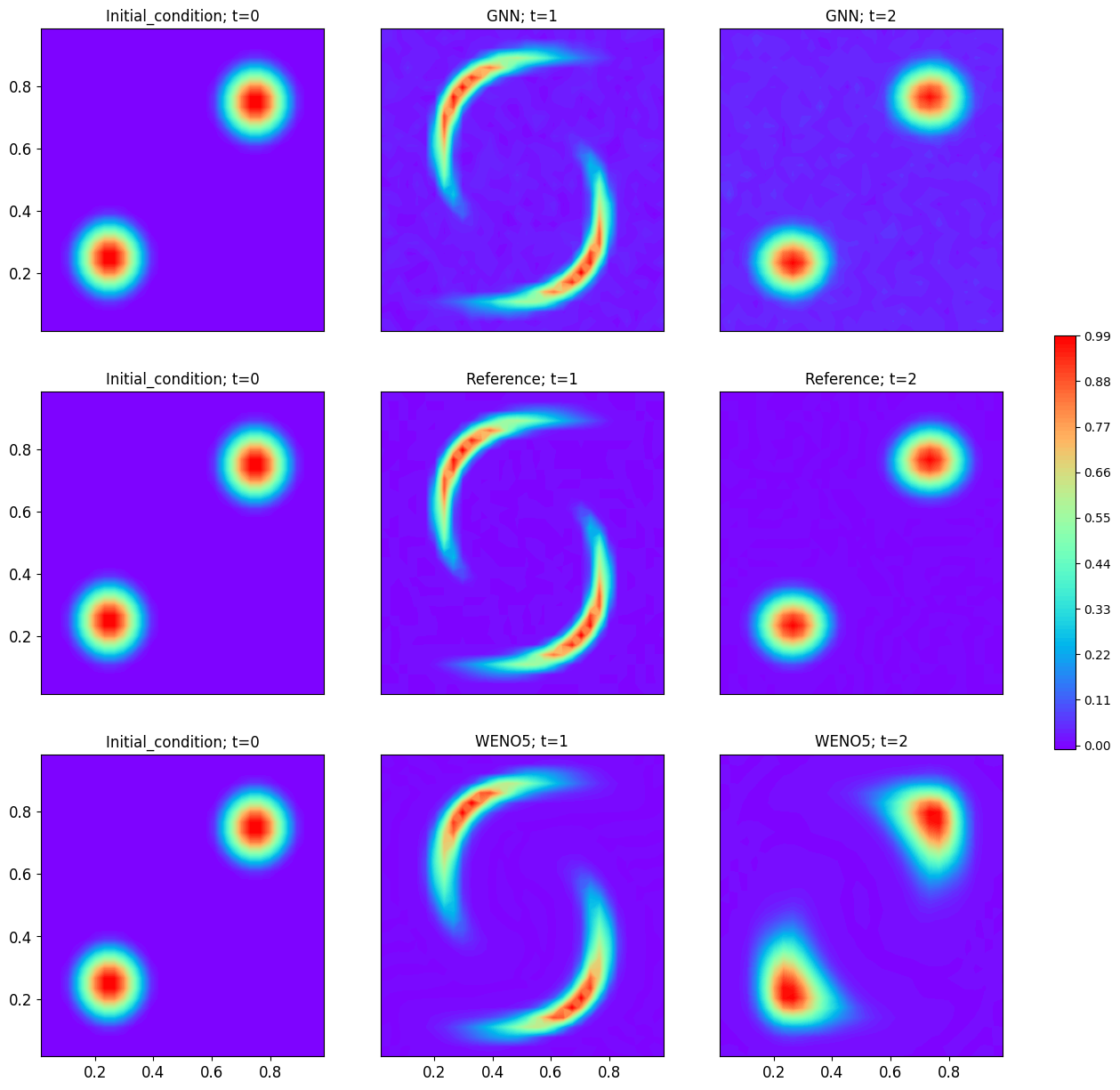}}
	\caption{Contour plots of the numerical solutions of the 2D deformational flow at $t=0,\,1,\,2$  in Example \ref{defor} for two cosine bells. CFL=10.2 for our method, CFL=0.6 for WENO5. $r_0=6,c_{1,x}=0.3,c_{1,y}=0.3,c_{2,x}=0.8,c_{2,y}=0.8$. The solver is trained with a data set for which each trajectory only contains a single cosine bell.}\label{fig:defor_onebell_test}
\end{figure}

Despite the training data consists of solution trajectories featuring a single bell, we demonstrate that the trained model can generalize to simulate problems with an initial condition that contained two randomly placed bells centered at [$c_{1,x},c_{1,y}$] and  [$c_{2,x},c_{2,y}$]:
\begin{equation}\label{eq:deformation-ini-two}
\begin{aligned}
    u(x,y) &= \frac{1}{2}[1+\cos (\pi r_1)+\cos (\pi r_2)]\\
    r_1(x,y) &= \min\left[1,r_0\sqrt{(x-c_{1,x})^2+(y-c_{1,y})^2}\right]\\
    r_2(x,y) &= \min\left[1,r_0\sqrt{(x-c_{2,x})^2+(y-c_{2,y})^2}\right].
\end{aligned} 
\end{equation}
In Figure \ref{fig:defor_onebell_test}, we present the contour plots of numerical solutions of the initial condition \eqref{eq:deformation-ini-two} with $c_{1,x}=0.3,c_{1,y}=0.3,c_{2,x}=0.8,c_{2,y}=0.8,r_0=6$. The CFL numbers are chosen as 10.2 and 0.6 for our method and WENO5, respectively. The reference solution is generated in the same way as the training data. Note that the solver can still capture deformations and accurately restore the initial profile at $t=T$.  In contrast, the solution computed using WENO5 exhibits severe smearing.

Lastly, we show the efficiency of the proposed method by providing the comparison of the errors and run-time between the proposed GNN-based model and the WENO5 method.  In Table
 \ref{table:defor-time}, we provide the run-time for simulating three test samples up to $t=T$ using our model with a mesh resolution of $32\times32$ cells and CFL=10.2, as well as the WENO5 with mesh resolutions of  $32\times32$ cells  and   $128\times128$ cells, both utilizing a CFL number of 0.6. Table \ref{table:defor-error} reports the corresponding mean square errors. The computational time of the proposed GNN-based model is slightly smaller than the RK WENO5 method with the same grid resolution $32\times 32$ cells and much smaller than the WENO5 method with a resolution of $128\times 128$ cells. Meanwhile, our method achieves significantly smaller errors compared to the  WENO5 method with the same mesh size, and only slightly larger errors than the WENO5 method using a mesh with much higher resolution.

\begin{table}[!htbp]
		\centering
		\caption{Run-time comparison of Example \ref{defor}. Run-time (seconds) measure for one period from $t=0$ to $t=T$ on a single GeForce RTX 3090 Ti GPU for the proposed GNN-based model, and on a CPU for the traditional Eulerian RK WENO5 method using Python.} \label{table:defor-time}
		% \smallskip
  \small
		\begin{tabular}{@{}cp{3.5cm}p{3.2cm}p{2cm}@{}}
			\hline
			 Samples\textbackslash Method &  WENO5($128\times 128$) &  WENO5($32\times 32$) & GNN($32\times 32$) \\
			\hline
		Sample = 0	    &  76.1154 & 4.7250 & 1.0628   \\  \hline
            Sample = 1 &  79.1125 &5.1799&1.0531 \\
			\hline
            Sample = 2 &  78.2025 &4.9729 &1.0647 \\
			\hline

		\end{tabular}
	\end{table}

\begin{table}[!htbp]
		\centering
		\caption{Mean square errors of Example \ref{defor}.} \label{table:defor-error}
		\small
		\begin{tabular}{@{}cp{3.5cm}p{3.2cm}p{2cm}@{}}
			\hline
		 Samples\textbackslash Method & WENO5($128\times 128$) &  WENO5($32\times 32$) & GNN($32\times 32$) \\
			\hline
           Sample = 0 & 3.1357E-6 &1.5792E-3 & 2.7746E-5 \\
			\hline
             Sample = 1 & 3.2170E-6 &1.7188E-3 & 2.4379E-5\\
			\hline
            Sample = 2 & 2.9835E-6 &1.6179E-3 & 2.8598E-5 \\
			\hline

		\end{tabular}
	\end{table}


% \Cref{fig:testfig} shows some example results. Additional results are
% available in the supplement in \cref{tab:foo}.

% \begin{figure}[htbp]
%   \centering
%   \label{fig:a}\includegraphics{lexample_fig1}
%   \caption{Example figure using external image files.}
%   \label{fig:testfig}
% \end{figure}

% \lipsum[51]

% \section{Discussion of \texorpdfstring{{\boldmath$Z=X \cup Y$}}{Z = X union Y}}

% \lipsum[76]

\subsection{Nonlinear Vlasov-Possion System}
In this subsection, we present the numerical results for simulating the nonlinear 1D1V VP system. We demonstrate the efficiency and accuracy of the proposed GNN-based SL FD method coupled with the second-order RKEI by comparing it to the traditional SL FV WENO5 scheme.  The training data and reference solutions are generated by the fourth-order conservative SL FV WENO scheme. 
%For the nonlinear 1D1V VP system, we choose the node feature for node $I_{ij}$ as $x_{ij} = (f_{ij},\xi_{ij},\eta_{ij}) \in \mathbb{R}^3$. 

\begin{example}\label{landau-damping}
In this example, we consider the Landau damping with the initial condition 
    \begin{equation}\label{eq:landau-damping}
        f(x,v,t=0)=\frac{1}{\sqrt{2\pi}}(1+\alpha \cos(kx))\exp\left(-\frac{v^2}{2}\right), \quad x\in[0,L], \quad v\in[-V_c,V_c],
    \end{equation}
    where $k=0.5$, $L=4\pi$, and $V_c=2\pi$. 
\end{example}

We generate the training data by coarsening 6 solution trajectories with a $256 \times 512$-cell grid by a factor of 8 in each dimension. The initial conditions for the solution trajectories are determined using \eqref{eq:landau-damping}, with $\alpha$ randomly selected from a uniform distribution in the range of $[0.05,0.45]$. Each solution trajectory contains a sequence of time steps from $t=0$ to $t=40$, with CFL=10.8. The reference solution is generated with the same approach used to create the training data.

\begin{figure}[!htbp]
\centering
    \centerline{\includegraphics[width=0.9\textwidth]{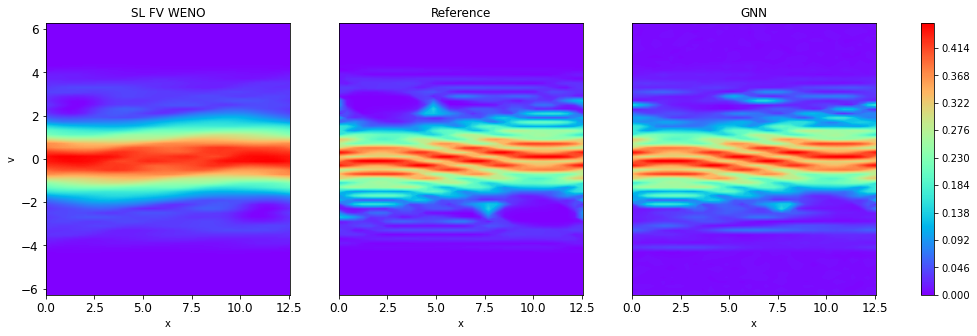}}
	\caption{Contour plots of numerical solutions of the strong Landau damping at $t=40$ in Example \ref{landau-damping} with $\alpha=0.5$.}\label{fig:landau_contour}
\end{figure}

\begin{figure}[!htbp]
    \centerline{\includegraphics[width=0.5\textwidth]{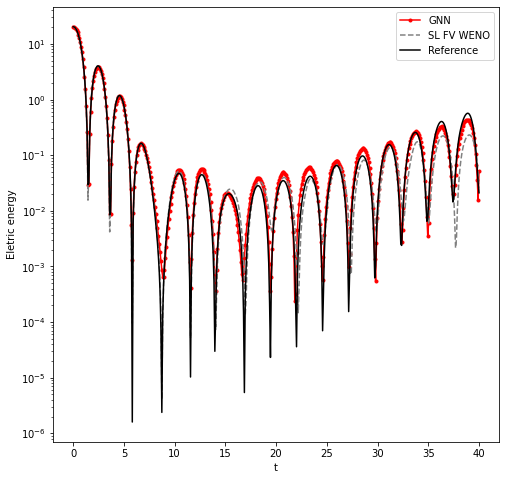}}
	\caption{Time histories of the electric energy of the strong Landau damping in Example \ref{landau-damping} with $\alpha=0.5$.}\label{fig:landau_eletric}
\end{figure}

After training, we test the performance of our solver for simulating the VP system with initial condition \eqref{eq:landau-damping} of $\alpha=0.5$, yielding the strong Landau damping, which lies outside the range of the training data. Figure \ref{fig:landau_contour} presents contour plots of the numerical solutions computed by our method and the traditional SL FV WENO method, all implemented on the same mesh with $32\times 64$ cells. It can be observed that our method can accurately capture the filamentation structure, and the results are in good agreement with the reference solution. The traditional SL FV WENO scheme produces reasonable results but exhibits smeared solution structures, mainly due to the low mesh resolution. We also plot the time histories of the electric energy for each approach in Figure \ref{fig:landau_eletric}. Our solver yields results that agree well with the reference solution, even better than the results of SL FV WENO with the same resolution.


\begin{example}\label{two-stream}
    In this example, we simulate the symmetric two stream instability with the initial condition 
    \begin{equation}\label{eq:two-stream}
        f(x,v,t=0)=\frac{1}{\sqrt{2\pi}}(1+\alpha \cos(kx))v^2\exp\left(-\frac{v^2}{2}\right),\quad  x\in[0,L],\quad v\in[-V_c,V_c],
    \end{equation}
    where $k=0.5$, $L=4\pi$, and $V_c=2\pi$. 
\end{example}

The training data is generated by coarsening 5 solution trajectories with a $256\times 512$ grid by a factor of 8 in each dimension. The initial conditions are determined using \eqref{eq:two-stream}, with $
\alpha$ randomly sampled from a uniform distribution in the range $[0.01,0.05]$. Each solution trajectory contains a sequence of time steps from $t=0$ to $t=53$. For the purpose of comparison, the reference is generated with the same approach used to create the training data.

\begin{figure}[!htbp]
\centering
    \centerline{\includegraphics[width=0.9\textwidth]{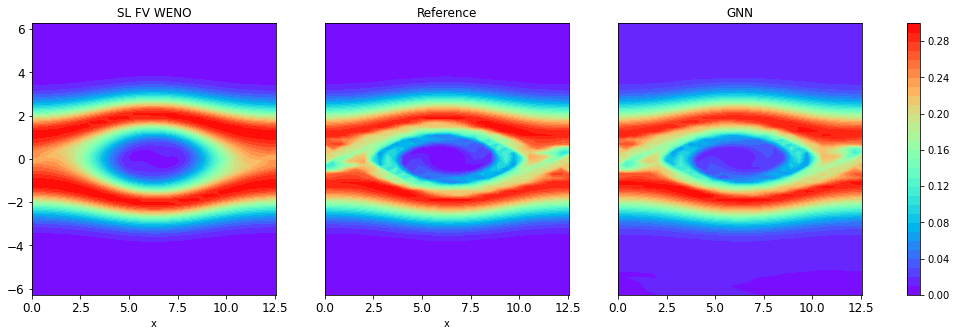}}
	\caption{Contour plots of numerical solutions of the two stream instability at $t=53$ in Example \ref{two-stream} with $\alpha=0.01$.  CFL = 10.8.}\label{fig:two_contour}
\end{figure}

\begin{figure}[!htbp]
    \centerline{\includegraphics[width=0.8\textwidth]{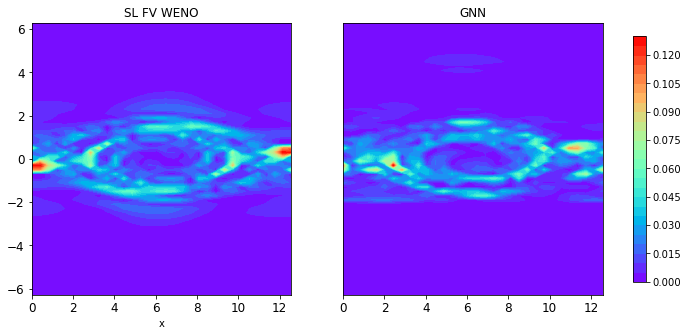}}
	\caption{The error of numerical solutions of the two stream instability at $t=53$ in Example \ref{two-stream} with $\alpha=0.01$. CFL = 10.8.}\label{fig:two_error}
\end{figure}

To test, we present the contour plots of the two-stream instability with $\alpha=0.01$ at $T=53$ in Figure \ref{fig:two_contour} for our method as well as the SL FV WENO scheme for comparison. Our method produces numerical results that are in good agreement with the reference solution, while the SL FV WENO scheme fails to capture fine-scale structures of interest, such as the roll-up at the center of the solution. We also plot the absolute error between the numerical solutions and the reference solution in Figure \ref{fig:two_error} to facilitate a more detailed comparison.

Furthermore, we demonstrate the efficiency of the proposed method by providing the comparison of  the run-time between the proposed GNN-based model and the  SL FV WENO method.  In Table
 \ref{table:two-time}, we provide the run-time for simulating three test samples up to $t=53$ using the proposed model with a mesh resolution of $32\times64$ cells and CFL number of 10.8, as well as the SL FV WENO with mesh resolutions of  $32\times64$ cells  and   $128\times256$ cells, both using the CFL number of 10.8.  
 % Table \ref{table:two-error} reports the corresponding mean square errors. 
It is observed that the computational time of the proposed multi-fidelity model is higher than that of the SL FV WENO method with the same mesh resolution of $32\times 64$ cells. However, it is lower than the SL FV WENO over a finer mesh of $128\times 256$ cells. We also note that implementing the SL FV WENO scheme requires substantial  human efforts, whereas our method benefits from a simpler code structure, thanks to the highly efficient ML packages utilized, including Pytorch \cite{paszke2019pytorch} and Pyg \cite{fey2019fast}.

% Meanwhile, the proposed model achieves much smaller errors compared to the SL FV WENO method with the same mesh size, and only slightly larger errors than the SLFD WENO method using the high-resolution mesh. 

\begin{table}[!htbp]
		\centering
		\caption{Run-time comparison of Example \ref{two-stream}. Run-time (seconds) measure for one period from $t=0$ to $t=53$ on a single GeForceRTX 3090 Ti GPU for the proposed GNN-based model, and on a CPU for the SL FV WENO method using Fortran.} \label{table:two-time}
		% \smallskip
  \small
		\begin{tabular}{@{}cp{3.5cm}p{3.2cm}p{2cm}@{}}
			\hline
			 Samples\textbackslash Method & SL FV WENO($128\times 256$) &  SL FV WENO($32\times 64$) & GNN($32\times 64$) \\
			\hline
		Sample = 0	    &  194.0781 & 2.7969 & 16.9319   \\  \hline
            Sample = 1 &  194.3593 &2.9844&16.6703 \\
			\hline
            Sample = 2 &  193.5625 &3.0000 &16.6089 \\
			\hline

		\end{tabular}
	\end{table}

% \begin{table}[!htbp]
% 		\centering
% 		\caption{Mean square errors of Example \ref{two-stream}.} \label{table:two-error}
% 		\small
% 		\begin{tabular}{@{}cp{3.5cm}p{3.2cm}p{2cm}@{}}
% 			\hline
% 		 Samples\textbackslash Method & SLFD WENO($128\times 256$) &  SLFD WENO($32\times 64$) & GNN($32\times 64$) \\
% 			\hline
%            Sample = 0 & 3.2369E-5 &2.5528E-4 & 3.9360E-5 \\
% 			\hline
%              Sample = 1 & 3.9071E-5 &2.1098E-4 & 3.5472E-5\\
% 			\hline
%             Sample = 2 & 3.8929E-5 &1.9640E-4 & 4.5432E-5 \\
% 			\hline

% 		\end{tabular}
% 	\end{table}


\begin{example}\label{two-stream-another}
In this example,  we consider another two stream instability with the following initial condition
    % \begin{equation}\label{eq:two-stream-another}
    %     f(x,v,t=0)=\frac{2}{7\sqrt{2\pi}}(1+5v^2)(1+\alpha_1\cos(kx)+\alpha_2\cos(2kx)+\alpha_3\cos(3kx))\mathrm{exp}(\frac{-v^2}{2}),  x\in[0,L], v\in[-V_c,V_c],
    % \end{equation}
    \begin{equation}
\begin{aligned}\label{eq:two-stream-another} &f(x, v, t=0) \\
& =  \frac{2}{7 \sqrt{2 \pi}}(1+5 v^2)(1+\alpha_1\cos(kx)+\alpha_2\cos(2kx)+\alpha_3\cos(3kx)) \exp \left(-\frac{v^2}{2}\right),\\
&x\in[0,L], v\in[-V_c,V_c],
\end{aligned}
\end{equation}
    where $k=0.5$, $L=4\pi$ and $V_c=2\pi$. This example is more challenging than the previous two examples, as the initial conditions are defined by perturbing the first three Fourier modes of the equilibrium, with amplitudes $\alpha_1$, $\alpha_2$, and $\alpha_3$, respectively.
\end{example}

We generate the training data by obtaining four solution trajectories on a $256\times 512$ grid and then downsampling them by a factor of 8 in each dimension. The initial conditions for each solution trajectory are given by \eqref{eq:two-stream-another}, with $\alpha_1, \alpha_2$ and $\alpha_3$ randomly sampled from a uniform distribution in the range $[0.01,0.02]$. Each solution trajectory contains a sequence of time steps from $t=0$ to $t=53$ with CFL=10.8. The reference solution is generated with the same approach used to create the training data. 

\begin{figure}[!htbp]
\centering
    \centerline{\includegraphics[width=0.9\textwidth]{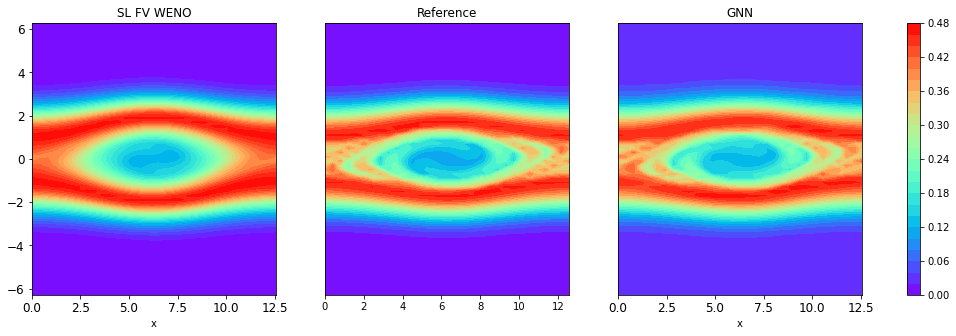}}
	\caption{Contour plots of numerical solutions of the two stream instability at $t=53$ in Example \ref{two-stream-another} with $\alpha_1=0.01, \alpha_2=0.01/1.2, \alpha_3=0.01/1.2$.  CFL = 10.8.}\label{fig:two1_contour}
\end{figure}

\begin{figure}[!htbp]
    \centerline{\includegraphics[width=0.8\textwidth]{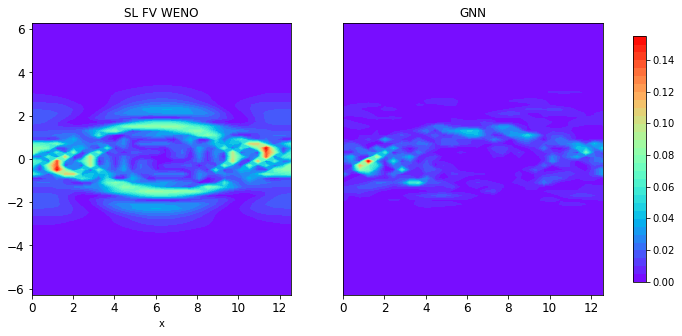}}
	\caption{The error of numerical solutions of the two stream instability at $t=53$ in Example \ref{two-stream-another} with $\alpha_1=0.01, \alpha_2=0.01/1.2, \alpha_3=0.01/1.2$. CFL = 10.8.}\label{fig:two1_error}
\end{figure}

During testing, we consider the initial condition \eqref{eq:two-stream-another} with $\alpha_1=0.01, \alpha_2=0.01/1.2$, and $\alpha_3=0.01/1.2$, which is a widely used benchmark configuration in the literature \cite{filbet2003comparison,qiu2011positivity,xiong2019conservative}. Note that such a parameter choice is outside the range of the training data. Figure \ref{fig:two1_contour} shows the contour plots of the numerical solutions computed by our method and the SL FV WENO scheme. It can be observed that the results by our solver qualitatively agree with the reference solution, effectively capturing
the underlying fine-scale structures of interest. Meanwhile, the SL FV WENO scheme can provide reasonable results, but tends to smear out the small-scale structures in the solution. This demonstrates that the proposed model can produce results with reasonable accuracy and possesses certain generalization capabilities. In Figure \ref{fig:two1_error}, we plot the difference between the numerical solutions by our method as well as by the SL FV WENO scheme and the reference solution. It can be observed that our method achieves smaller errors.


\section{Conclusions}
\label{sec:conclusions}

Semi-Lagrangian (SL) schemes are known as efficient numerical tools for simulating transport equations and have been widely used across various fields. Despite the effectiveness of the SL methodology, directly designing conservative SL schemes for simulating multidimensional problems remains a significant challenge. In this paper, we proposed a novel conservative ML-based SL finite difference scheme that allows for extra-large time step evolution. Our method employs an end-to-end neural architecture  to learn the optimal SL discretization through message passing of the GNN. The core of the solver is to construct a dynamical graph that can handle the complexities associated with accurately tracking upstream points along characteristics. By learning a conservative SL discretization via a data-driven approach, our method can achieve improved accuracy and efficiency compared to traditional numerical schemes. Meanwhile, it can simplify the implementation of SL algorithms. Numerical examples conducted in this work, including benchmark transport equations in both one and two dimensions and nonlinear Vlasov-Poisson system, demonstrate the efficiency of the proposed method. Future work includes improving the generalization capabilities of the method, addressing the challenges related to adaptivity and complex geometries, and exploring the possibility of applying the method to other problems, including convection dominated equations.

% \appendix
% \section{An example appendix} 
% % \lipsum[71]

% \begin{lemma}
% Test Lemma.
% \end{lemma}

% \section*{Acknowledgments}
% We would like to acknowledge the assistance of volunteers in putting
% together this example manuscript and supplement.

\bibliographystyle{siamplain}
\bibliography{references}
\end{document}

\maketitle

\section{A detailed example}

Here we include some equations and theorem-like environments to show
how these are labeled in a supplement and can be referenced from the
main text.
Consider the following equation:
\begin{equation}
  \label{eq:suppa}
  a^2 + b^2 = c^2.
\end{equation}
You can also reference equations such as \cref{eq:matrices,eq:bb} 
from the main article in this supplement.

\lipsum[100-101]

\begin{theorem}
  An example theorem.
\end{theorem}

\lipsum[102]
 
\begin{lemma}
  An example lemma.
\end{lemma}

\lipsum[103-105]

Here is an example citation: \cite{KoMa14}.

\section[Proof of Thm]{Proof of \cref{thm:bigthm}}
\label{sec:proof}
\lipsum[106-112]

\section{Additional experimental results}
\Cref{tab:foo} shows additional
supporting evidence. 

\begin{table}[htbp]
{\footnotesize
  \caption{Example table}  \label{tab:foo}
\begin{center}
  \begin{tabular}{|c|c|c|} \hline
   Species & \bf Mean & \bf Std.~Dev. \\ \hline
    1 & 3.4 & 1.2 \\
    2 & 5.4 & 0.6 \\ \hline
  \end{tabular}
\end{center}
}
\end{table}

\bibliographystyle{siamplain}
\bibliography{references}